\documentclass[proc]{edpsmath}

\usepackage{graphicx}
\usepackage{placeins} 

\begin{document}

\title{Hastings-Metropolis algorithm on Markov chains for small-probability estimation}\thanks{The work presented in this manuscript was carried out in the framework of the REDVAR project of the CEMRACS 2013}\thanks{This work was financed by the Commissariat \`a l'energie atomique et aux \'energies alternatives}
\author{F. Bachoc}\address{Department of Statistics and Operations Research, University of Vienna, Oskar-Morgenstern-Platz 1, A-1090 Vienna}
\secondaddress{When the work presented in this manuscript was carried out, the author was affiliated to CEA-Saclay, DEN, DM2S, STMF, LGLS, F-91191 Gif-Sur-Yvette, France  and to the Laboratoire de Probabilit\'es et Mod\`eles Al\'eatoires, Universit\'e Paris VII}
\author{A. Bachouch}\address{Laboratoire Manceau de Math\'ematiques (LMM), Universit\'e du Maine. Avenue Olivier Messiaen 72085 Le Mans CEDEX 9 . France}
\secondaddress{Laboratoire de Mod\'elisation Math\'ematique et num\'erique dans les sciences de l'ing\'enieur (Lamsin). ENIT. BP 37, 1002 Tunis. Tunisie}
\author{L. Len\^otre}\address{Inria, Research Centre Rennes-Bretagne Atlantique, Campus de Beaulieu 35042 Rennes Cedex France} 
\secondaddress{Universit\'e de Rennes 1, Campus de Beaulieu 35042 Rennes Cedex France}
%
%
\begin{abstract} Shielding studies in neutron transport, with Monte Carlo codes, yield challenging problems of small-probability estimation. The particularity of these studies is that the small probability to estimate is formulated in terms of the distribution of a Markov chain, instead of that of a random vector in more classical cases. Thus, it is not straightforward to adapt classical statistical methods, for estimating small probabilities involving random vectors, to these neutron-transport problems. A recent interacting-particle method for small-probability estimation, relying on the Hastings-Metropolis algorithm, is presented. It is shown how to adapt the Hastings-Metropolis algorithm when dealing with Markov chains. A convergence result is also shown.
Then, the practical implementation of the resulting method for small-probability estimation is treated in details, for a Monte Carlo shielding study. Finally, it is shown, for this study, that the proposed interacting-particle method considerably outperforms a simple-Monte Carlo method, when the probability to estimate is small.
 \end{abstract}
\begin{resume} Dans les \'etudes de protection en neutronique, celles fond\'ees sur des codes Monte-Carlo posent d'importants probl\`emes d'estimation de faibles probabilit\'es. La particularit\'e de ces \'etudes est que les faibles probabilit\'es sont exprim\'ees en termes de lois sur des chaines de Markov, contrairement \`a des lois sur des vecteurs al\'eatoires dans les cas les plus classiques. Ainsi, les m\'ethodes classiques d'estimation de faibles probabilit\'es, portant sur des vecteurs al\'eatoires, ne peuvent s'utiliser telles qu'elles, pour ces probl\`emes neutroniques. Un m\'ethode r\'ecente d'estimation de faibles probabilit\'es, par syst\`eme de particules en int\'eraction, reposant sur l'algorithme de Hastings-Metropolis, est pr\'esent\'ee. Il est alors montr\'e comment adapter l'algorithme de Hastings-Metropolis au cas des chaines de Markov. Un r\'esultat de convergence est ainsi prouv\'e. Ensuite, il est expliqu\'e en d\'etail comment appliquer la m\'ethode obtenue \`a une \'etude de protection par Monte-Carlo. Finalement, pour cette \'etude, il est montr\'e que la m\'ethode par syst\`eme de particules en int\'eraction est consid\'erablement plus efficace qu'une m\'ethode par Monte Carlo classique, lorsque la probabilit\'e \`a estimer est faible. \end{resume}
\maketitle
\section{Introduction} \label{section:intro}

The study of neutronics began in the 40’s, when nuclear energy was on the verge of being used both for setting up nuclear devices like bombs and for civil purposes like the production of energy.
Neutronics is the study of neutron population in fissile media that can be modeled using the linear Boltzmann equation, also known as the transport equation.
More precisely, it can be subdivided in two different sub-domains. On the one hand, criticality studies aim at understanding the neutron population dynamics due to the branching process 
that mimics fission reaction (see for instance \cite{zoaone} for a recent survey on branching processes in neutronics). On the other hand, when neutrons are propagated through media where fission reactions do not occur,
or can safely be neglected, their transport can be modeled by simple exponential flights \cite{zoatwo}: indeed, between each collisions, neutrons travel along straight path distributed exponentially. 

Among this last category, shielding studies allow to size shielding structures so as to protect humans from ionizing particles, and imply, by definition, the attenuation of initial neutron
flux typically by several decades. For instance, the vessel structure of a nuclear reactor core attenuates the thermal neutron flux inside the core by a factor roughly equal to $10^{13}$.
Many different national nuclear authorities require shielding studies of nuclear systems before giving their agreement for the design of these systems. Examples are reactor cores, but also devices for nuclear medicine (proton-therapy, gamma-therapy, etc…). The study of those nuclear systems is complicated by 3-dimensional effects due to the geometry and by non-trivial energetic spectrum that can hardly be modeled.

Since Monte Carlo transport codes (like MCNP \cite{mcnp}, Geant4 \cite{ago}, Tripoli-4 \cite{diop}) require very few hypotheses,
they are often used for shielding studies. 
Nevertheless, those studies represent a long-standing numerical challenge for Monte Carlo codes in the sense that they schematically require to evaluate the proportion of neutrons that ``pass'' through the shielding system. This proportion is, by construction, very small.
Hence a shielding study by Monte Carlo code requires to evaluate a small probability, which is the motivation of the present paper.

There is a fair amount of literature on classical techniques for reducing the variance in these small-probability estimation problems for Monte Carlo codes.
Those techniques often rely on a zero-variance scheme \cite{kahn,hoog,booone} adapted to the Boltzmann equation, allied with “weight-watching” techniques \cite{bootwo}.
The particular forms that this scheme takes when concretely developed in various transport codes range from the use of weight windows \cite{burn,hoog,kahn,mcnp}, like in MCNP, to the use of the exponential transform \cite{both90automated,diop} like in Tripoli-4.  
Nowadays, all those techniques have proven to be often limited in view of fullfilling the requirements made by national nuclear authorities for the precise measurements of radiation, which standards are progressively strengthened. 
Thus, new variance reduction techniques have been recently proposed in the literature (see for instance \cite{dum} for the use of neural networks for evaluating the importance function).

This paper deals with the application of the recent interacting-particle method developed in \cite{guyader11simulation}, which has interesting theoretical properties and is particularly efficient in various practical cases. 
Nevertheless, the application is not straightforward, since the method in \cite{guyader11simulation} is designed for finite-dimensional problems, while the output of a neutron transport Monte Carlo code consists in a trajectory of the stochastic process driving the neutron behavior.

More specifically, as a neutron travels along straight paths between collisions, 
there is no loss of information in considering only the characteristics of the collisions (dates, positions, energies, subparticle creations) as random. Furthermore, in order to simplify the matter, we choose to only consider the simple but realistic case of a monokinetic particle (constant speed) and to avoid the problem of the subparticle-creation phenomena and the energy dependence.
As a result, the dynamics of the particle consists of a Markov chain, whose paths are the sets of successive collisions, and such that it is absorbed after a finite number of collisions. 
This model already implies serious complications for the interacting-particle method which relies on the Hastings-Metropolis algorithm \cite{metro,hast} for practical implementation. 
Indeed, this is not automatic to apply this algorithm to sample the paths of a Markov chain that is absorbed almost surely instead of just simulating a random vector.
The main difficulty relies on the randomness of the path length which compels the algorithm to jump through finite dimensional spaces.
An algorithm was proposed for a problem with also dimensional jumps which occurs in Bayesian model determination \cite{grn}.
However, the law to sample was simpler than the global distribution of a Markov chain and the state space totally different. 
Thus, our contribution is two-fold, as first we show how the Hastings-Metropolis algorithm can be extended to the case of Markov chains that are absorbed after finite time, and second we adapt the resulting interacting-particle  method to the Monte Carlo simulation of a monokinetic particle in a simplified but realistic shielding system.
We perform several numerical simulations which show that the smaller the probability to estimate is, the more the method we propose clearly outperforms a simple Monte Carlo method.

In what follows, we start with a short introduction to the interacting-particle method \cite{guyader11simulation} and highlight the need of the Hastings-Metropolis algorithm (Section \ref{section:reminder:AG}). 
Then, we dedicate a consequent work to prove the validity and the convergence of the Hastings-Metropolis algorithm applied to Markov Chains very similar to the one used in neutronic in order to convince the reader of the adaptation (Section \ref{section:HM:on:Markov:chains}).
After that, we present the aforementioned monokinetic particle model and give the actual equations for the small probability estimation method
\footnote{The reader interested in the neutronic Monte Carlo application can go directly from Section \ref{section:reminder:AG} to Section \ref{section:one:two:dim:case}.} (Section \ref{section:one:two:dim:case}). 
At last, we show the obtained numerical results for shielding studies and discuss them (Section \ref{section:numerical:results}).

\section{The Interacting-Particle Method for Small Probability Estimation} \label{section:reminder:AG}

Let $(\Omega,\mathcal{F},P)$ be probability space, $(S,\mathcal{S},Q)$ a measured space, $X$ a random variable from $(\Omega,\mathcal{F},P)$ to $(S,\mathcal{S},Q)$ that can be sampled and $\Phi: S \to \xR$ an objective function with a continuous cumulative distribution function $F$. 
We aim at estimating the probability $p$ of the event $\Phi(X) \geq l$, for a given level $l \in \xR$.
In order to evaluate $p$, we choose to use the interacting-particle method introduced in \cite{guyader11simulation}. 

\subsection{Theoretical Version of the Interacting-Particle Method} \label{subsection:ag:ideal}

Let assume that we are able to sample $X$ conditionally to the event $\Phi(X) \geq x$, for any $x \in \xR$.  
In this case, the interacting-particle method \cite{guyader11simulation} for estimating $p$, described in the following algorithm, yields a conceptually convenient estimator $\hat{p}$ with explicit finite-sample distribution.
The algorithm is parameterized by the number of particles $N$.

\begin{lgrthm} \label{algo:AG}

~

\begin{itemize}
\item Generate an $iid$ sample $(X_1,..,X_N)$, from the distribution of $X$, and initialize $m=1$, \\
 $L_1= \min(\Phi(X_1),..,\Phi(X_N))$ and $X^1_1=X_1,...,X^1_N=X_N$.
\item While $L_m \leq l$ do
\begin{itemize}
\item For $i=1,...,N$
\begin{itemize}
\item Set $X^{m+1}_i = X^{m}_i$ if $\Phi(X^{m}_i)>L_m$, and else $X^{m+1}_i = X^*$, where $X^*$ follows the distribution of $X$ conditionally to $\Phi(X) \geq L_m$, and is independent of any other random variables involved in the algorithm.
\end{itemize}
\item Set $m=m+1$.
\item Set $L_m=\min(\Phi(X^m_1),..,\Phi(X^m_N))$.
\end{itemize}
\item The estimate of the probability $p$ is $\hat{p}=(1-\frac{1}{N})^{m-1}$.
\end{itemize}

\end{lgrthm}

For each $N < \infty$, the estimator $\hat{p}$ has an explicit distribution which is detailed in \cite{guyader11simulation}. 
This reference exhibits two properties of $\hat{p}$ : the estimator is unbiased and an asymptotic $95 \%$ confidence interval, for $N$ large, has the following form
\begin{equation} \label{eq:IC95:AG}
I_{\hat{p}} = 
\left[
 \hat{p} \exp \left( -1.96 \sqrt{  \frac{-\log{\hat{p}}}{N} } \right)
,
\hat{p} \exp \left( 1.96 \sqrt{  \frac{-\log{\hat{p}}}{N} } \right)
  \right].  
\end{equation}
We note that the event $p \in I_{\hat{p}}$ is asymptotically equivalent to the event
$\hat{p} \in I_{p}$, with $I_p$ as in \eqref{eq:IC95:AG} with $\hat{p}$ replaced by $p$. We mean that the probabilities of the two events converge to $0.95$ and the probability of their symmetric difference converges to $0$ as $N \to \infty$. The asymptotic equivalence holds because $\log(\hat{p})$ is asymptotically normal with mean $\log(p)$ and variance $-\log(p)/N$ \cite{guyader11simulation}.
We will use this property in Section \ref{section:numerical:results}.

\subsection{Practical implementation using the Hastings-Metropolis algorithm} \label{subsection:ag:practical}

In many practical cases, the previous algorithm is inapplicable, as it relies on the strong assumption of being able to exactly sample $X$ conditionally to $\Phi(X) \geq t$, for any $t \in \xR$.
Subsequently, the authors of \cite{guyader11simulation} propose to use the Hastings-Metropolis algorithm to simulate this conditional distribution. This method requires the following assumptions:

\begin{itemize}
\item The distribution of $X$ has a probability distribution function (pdf) $f$ with respect to
$(S,\mathcal{S},Q)$. For any $x \in S$ we can compute $f(x)$.
\item We have available a transition kernel on $(S,\mathcal{S},Q)$ with conditional pdf $\kappa(x,y)$ (pdf at $y$ conditionally to $x$). [Throughout this paper, $\kappa$ is called the instrumental kernel.]
 We are able to sample from $\kappa(x,.)$ for any $x \in S$ and we can compute $\kappa(x,y)$ for any $x,y \in S$.
\end{itemize}

Let $t \in \xR$ and $x \in S$ so that $\Phi(x) \geq t$. Then, the following algorithm enables to, starting from $x$, sample approximately with the distribution of $X$, conditionally to $\Phi(X) \geq t$. The algorithm is parameterized by a number of iterations $T \in \xN^*$.
\begin{lgrthm} \label{algo:HM}

~

\begin{itemize}
\item Let $X = x$.
\item For $i=1,...,T$
\begin{itemize}
\item Independently from any other random variable, generate $X^*$ following the $\kappa(X,.)$ distribution.
\item If $\Phi(X^*) \geq t$
\begin{itemize}
\item Let $r = \frac{f(X^*)\kappa(X^*,X)}{f(X)\kappa(X,X^*)}$.
\item With probability $min(r,1)$, let $X=X^*$.
\end{itemize}
\end{itemize}
\item Return $X_{T,t}(x) = X$.
\end{itemize}

\end{lgrthm}

For consistency, we now give the actual interacting-particle method, involving Algorithm \ref{algo:HM}. This method is parameterized by the number of particles $N$ and the number of HM iterations $T$.

\begin{lgrthm} \label{algo:AG:with:HM}

~

\begin{itemize}
\item Generate an $iid$ sample $(X_1,..,X_N)$ from the distribution of $X$ and initialize
$m=1$, $L_1= \min(\Phi(X_1),..,\Phi(X_N))$ and $X^1_1=X_1,...,X^1_N=X_N$.
\item While $L_m \leq l$ do
\begin{itemize}
\item For $i=1,...,N$
\begin{itemize}
\item If $\Phi(X^{m}_i)>L_m$, set $X^{m+1}_i = X^{m}_i$.
\item Else sample an integer $J$ uniformly in the set $\{ 1 \leq j \leq N ; \Phi(X^{m}_j)>L_m \}$. Apply Algorithm \ref{algo:HM} with number of iterations $T$, starting point $X_J^m$ and with threshold value $t = L_m$. Write $X_{T,L_m}(X^m_J)$ for the output of this algorithm and let $X^{m+1}_i = X_{T,L_m}(X^m_J)$.
\end{itemize}
\item Set $m=m+1$.
\item Set $L_m=\min(\Phi(X^m_1),..,\Phi(X^m_N))$.
\end{itemize}
\item The estimate of the probability $p$ is   $\hat{p}=(1-\frac{1}{N})^{m-1}$.
\end{itemize}

\end{lgrthm}

The estimator $\hat{p}$ of Algorithm \ref{algo:AG:with:HM} is the practical estimator that we will study in the numerical results of Section \ref{section:numerical:results}. 
In \cite{guyader11simulation}, it is shown that, when the space $S$ is a subset of $\xR^d$, under mild assumptions, the distribution of the estimator of Algorithm \ref{algo:AG:with:HM} converges, as $T \to +\infty$, to the distribution of the ideal estimator of Algorithm \ref{algo:AG}. For this reason, we call the estimator of Algorithm \ref{algo:AG} the estimator corresponding to the case $T = +\infty$. We also call the confidence intervals \eqref{eq:IC95:AG} the confidence intervals of the case $T = + \infty$.

Nevertheless, as we discussed in Section \ref{section:intro}, the space $S$ we are interested in is a space of sequences that are killed after a finite time. Thus, it is not straightforward that
the convergence, as $T \to + \infty$, discussed above, hold in our case. Furthermore, even the notion of pdf on this space of sequences has to be defined.
This is the object of the next Section \ref{section:HM:on:Markov:chains}, that defines the notion of pdf, on a space of sequences that are killed after a finite time, and that gives a convergence result for the HM algorithm.
The definition of the pdf is also restated in Section \ref{section:one:two:dim:case}, so that Sections \ref{section:reminder:AG} and \ref{section:one:two:dim:case} are self-sufficient for the implementation of the small-probability estimation method for the monokinetic particle simulation.

\section{An extension of Hastings-Metropolis algorithm to Path Sampling} \label{section:HM:on:Markov:chains}

\subsection{Reformulation of the Markov Chain Describing the Neutronic Problem}

In many neutronic models, the dynamics of the collisions are described by a Markov chain $(X_n)_{n \geq 0}$ with values in $\mathbb{R}^d$ and which possesses a probability transition function $q$ and an initial position $X_0=\alpha$.
Since the detection problem occurs only in a restricted area, we decide to change this description using a censorship.
Such a trick will be of great help for the theoretical treatment developed later.

Let $D$ be an open bounded subset of $\xR^d$ with $\partial D$ its boundary. 
Because $D$ is the domain of interest, we rewrite the transition function of the process $(X_n)_{n \geq 0}$ as follows
 \begin{equation*}
 k(x,dy) = (q(x,y) 1_D(y) \, dy + q_x(D^C) \delta_{\Delta}(dy))1_{D}(x) + \delta_{\Delta}(dy) 1_{\Delta}(x)
\end{equation*}
where $\Delta$ is a resting point and 
\begin{equation*}
 q_x(D^C) = \int_{D^C} q(x,y) \, dy.
\end{equation*}
This kernel describes the following dynamic: 
\begin{itemize}
 \item while $(X_n)_{n \geq 0}$ is inside $D$, it behaves with the transition kernel q that reflects the dynamics of the collision and can push the neutron outside $D$ .
 \item when $(X_n)_{n \geq 0}$ enters in $D^C$, it is killed and sent to the resting point $\Delta$ where it stays indefinitely. This way we keep only the informationss occuring exactly inside $D$.
\end{itemize}
We call this stochastic process a Boundary Absorbed Markov Chains (BAMC).

\subsection{Reminder of the Hastings-Metropolis Algorithm}

The Hastings-Metropolis algorithm is a generic procedure used to sample a distribution $\gamma$ that admits a density with respect to a measure $\Pi$ \cite{metro,hast}.
The idea of this algorithm is to define a Markov chain $(Y_n)_{n \geq 0}$ with a transition kernel $\Gamma$ that converges to $\gamma$ in some sense that will be discussed later. 
In order to construct $(Y_n)_{n \geq 0}$, the Hastings-Metropolis procedure uses an instrumental Markov chain $(Z_n)_{n \geq 0}$ and an acceptation-rejection function $r$.
We will denote by $\kappa$ the probability transition kernel of $(Z_n)_{n \geq 0}$ and call it the instrumental kernel.
The main hypothesis required on $\kappa$ and $\Gamma$ fby the algorithm is that they admit a density with respect to the measure $\Pi$. 
A Step by step description of the algorithm: 
\begin{itemize}
 \item Introduce a starting point $x$ and use it to sample a potential new position $y$ of $(Z_n)_{n \geq 0}$.
 \item Accept $y$ and set $x=y$ or reject it using $r$.
 \item return the position $x$ as the sample.
\end{itemize}
The more this procedure is repeated the more approximation is reliable. We can write the transition kernel $\Gamma$ of $(Z_n)_{n \geq 0}$ as follows
\begin{equation*}
\Gamma(u,dv)=\overline{\kappa}(u,v)\Pi(dv)+\overline{r}(u)\delta_u(dv) 
\end{equation*}
where
\begin{equation*}
 \overline{\kappa}(u,v)=
 \begin{cases}
  \kappa(u,v) r(u,v),\quad &\text{if } x\not= y, \\
  0,&\text{if } x=y,
 \end{cases}
\end{equation*}
and
\begin{equation*}
 \overline{r}(u)=1-\int \overline{\kappa}(u,v) \Pi(dv).
\end{equation*}

We conclude this reminder of the Hastings-Metropolis algorithm with an example for the acceptation-rejection function $r$.
This function is of major importance as it impacts the efficiency of the procedure and ensures that $\Pi$ is invariant for $\Gamma$ \cite{tier}.
we refer to the literature for more details \cite{hast, metro, meyn, tier}.
Since it possesses a reversibility property that quickly provides the condition just mentionned above \cite{tier}, we propose: 
\begin{equation*}
 r(u,v)=
 \begin{cases}
  \min\left\lbrace \dfrac{\gamma(v) \kappa(v,u)}{\gamma(u) \kappa(u,v)},1 \right\rbrace, &\text{if }\gamma(u) \kappa(u,v)>0\\
  1, &\text{if }\gamma(u) \kappa(u,v)=0
 \end{cases}.
\end{equation*}

\subsection{Definition of a Point Absorbed Markov Chain on a Sphere} \label{subsection:distribution:definition:PAMC}

The extension of the Hastings-Metropolis to sample the paths of a BAMC is quite natural and has been already used in several numerical methods.
But, as far as we know, there is still no rigorous proof for the convergence.
As a result, we propose to conduct a proof for the below defined Point Absorbed Markov Chain (PAMC) on a sphere S that can be linked through differential geometry to the BAMC presented earlier .

Let $S$ be a sphere included in the unit ball of $\mathbb{R}^d$ and containing $0$.
We denote by $S_0$ the subset $S-\{0\}$ and by $\lambda$ the Lebesgue measure restrained on $S_0$.
Let remark that $\lambda$ is the same on both $S_0$ and $S$, and that the densities are identical.
As a result, $\lambda$ will also stand for the Lebesgue measure on $S$.
We define a PAMC on S as the stochastic process $(M_n)_{n \geq 0}$ with value in $S_0$ and a probability transition function $m$ of the form 
\begin{equation*}
 m(x,dy) = (p(x,y) 1_{S_0}(y) \, dy + P_x(0) \delta_{0}(dy))1_{S_0}(x) + \delta_{0}(dy) 1_{0}(x).
\end{equation*}
where
\begin{enumerate}
	\item $p$ is a transition function on $S_0$ having a density with respect to $\lambda$,
	\item $(P_x(0))_{x \in S_0}$ is a family of positive real numbers,
	\item For every $x \in S_0$,
		\begin{equation*}
 			\int_{S} p(x,y) 1_{S_0}(y) \, dy + P_x(0) \delta_{0}(dy) = 1,
		\end{equation*}
	\item $m$ is such that $(M_n)_{n \geq 0}$ is almost surely absorbed in finite time.
	\item For every $(x,y) \in S_0^2$,
		\begin{equation*}
 			m < p(x,y) < M,
		\end{equation*}
	\item For every $x\in S_0$,
		\begin{equation*}
 			P_x(0) > c.
		\end{equation*}
\end{enumerate}

The proof of the Hastings-Metropolis algorithm extended to the above PAMC will be performed using results provided by some classical references \cite{meyn, num, tier}. 
We recall that the proof consists of showing that the Markov chain constructed with the Hastings-Metropolis procedure converges with respect to the total variation norm to the law of the Markov chain $(M_n)_{n \geq 0}$.
In order to use these results which suppose that the state space of the Hastings-Metroplis chain is a locally compact and separable topological space equipped with a metric, we have to consider few topological questions.

We start by pointing out that the state space of a PAMC on $S$ is the space of sequences with values in $S$ that are almost zero and which will be denoted by:
\begin{equation*}
 c_0(S)=\{(u_n)_{n \geq 0} \in S^{\mathbb{N}}: \exists n_0 \in \mathbb{N}, \forall n \geq n_0, u_n=0\},
\end{equation*}
We equipped this space with the norm:
\begin{equation*}
\|u\|_{\infty}=\underset{n \geq 0}{\text{max}} \, \|u_n\|_{\mathbb{R}^d}.
\end{equation*}
This state space have the properties mentioned earlier if we accept the following:
\begin{clm}
	It exists a locally compact and separable topology on the space $c_0(S)$ that can be equipped with a metric.
	In addition, the Borel $\sigma$-algebra generated by this topology coincides with $\mathcal{B}(c_0(S))$ the one generated by $\| \cdot \|_{\infty}$.
\end{clm}
\begin{rmrk}
	In order to share the intuition we had when we formulated this claim, we provide few arguments to the reader.
	Firstly, $c_0(S)$ is a subspace of $l^{\infty}(\mathbb{R}^d)$ which is locally compact when we set the weak-star topology.
	Secondly, this topology restrained to the unit ball can be equipped with a metric.
	Finally, $c_0(S)$ can be continuously injected in the space $l^{2}(\mathbb{R}^d)$ and it has been proved that the Borel $\sigma$-algebra generated by the weak-star topology on a Hilbert space
	coincides with the one generated by the topology of the norm.
\end{rmrk}

We start by mentioning that a PAMC on $S$ is a random variable:
\begin{eqnarray*}
  M: (\Omega,\mathcal{F},\mathbb{P}) & \mapsto & (c_0,\mathcal{B}(c_0(S))) \\
  \omega & \mapsto & (M_n(\omega))_{n \geq 0},
\end{eqnarray*}
if we use the $\sigma$-algebra:
\begin{equation*}
 \mathcal{F}=\underset{i=0}{\overset{+\infty}{\bigotimes}} \mathcal{B}(S),
\end{equation*}
generated by the Borelian cylinders of finite dimension.
Therefore, the following result shows the measurability of the process $(M_n)_{n \geq 0}$ with respect to the Borel $\sigma$-algebra generated by $\|\cdot\|_{\infty}$.
\begin{prpstn}
	The trace $\sigma$-algebra $\mathcal{F}_{|c_0(S)}$ on the subspace $c_0(S)$ of $\mathcal{F}$ is equal to $\mathcal{B}(c_0(S))$.
\end{prpstn}
\begin{proof}
 	Let $p_n$ be the projection from $c_0(S)$ in $S$ which associates $u_n$ to $u$. 
 	This application is Lipschitz. In fact, let $u$ and $v$ be in $c_0(S)$,
 	we have $\|u_n-v_n\|_{\mathbb{R}^d} \leq \|u-v\|_{\infty}$. Consequently, every projection is
 	measurable and we have the following inclusion:
	\begin{equation*}
 		\mathcal{F}_{|c_0(S)} \subset \mathcal{B}(c_0(S)).
	\end{equation*}
	
	On the other side, we know that $\mathcal{B}(c_0(S))$ is generated by the balls of radius $\rho \in \mathbb{Q}$ and center points $u \in T$ where $T$ is a dense subset of $S$, since $c_0(S)$ 			equipped with the norm $\| \cdot \|_{\infty}$ is separable.
 	Thus, it is enough to show that the ball $B(\rho,u)$ is in $\mathcal{F}_{|c_0(S)}$. 
 	In order to prove that, we write:
 	\begin{equation*}
  		B(\rho,u)=\overset{+\infty}{\underset{n=0}{\bigcap}} \{v \in c_0(S), \|u_n-v_n\|_{\mathbb{R}^d} \leq \rho \}.
 		\end{equation*}
 	Since each member of this intersection is in $\overline{\mathcal{F}}_{c_0}$, we have the opposite inclusion:
 	\begin{equation*}
 		\mathcal{B}(c_0(S)) \subset \overline{\mathcal{F}}_{c_0(S)}.
	\end{equation*} 
\end{proof}
\begin{rmrk}
	This proof can be considered as an adaption of a classical result for the Brownian Motion \cite{cmts}.  
\end{rmrk}

\subsection{Density of a Point Absorbed Markov Chain on a Sphere}

In order to use the Hastings-Metropolis algorithm, we must show that the law of a PAMC on the sphere $S$ admits a density with respect to a measure on $c_0(S)$.
Since we deal with a Markov process, we do not have to take the initial law into account.
As a result, we just have to find a density for the law of the process conditioned to start from $a$.

Without loss of generality, we can shift the element of $c_0(S)$ and rewrite them $(u_n)_{n \geq 1}$. 
Let introduce a partition of the space $c_0(S)$ using the subsets $(A_n)_{n \geq 0}$ consisting of:
\begin{equation*}
 A_0=\{u \in c_0(S) : u_k = 0, \forall  k \geq 1\},
\end{equation*}
\begin{equation*}
 A_n=\{u \in c_0(S) : u_k \in S_0, \forall  k \leq n \text{ and } u_k=0, \forall k > n\}, \forall n \geq 1
\end{equation*}
and the family of applications $(\pi_n)_{n \geq 0}$ defined as:
 \begin{equation*}
 \begin{split}
  \pi_n: c_0(S) &\mapsto S^n \\
  (u_n)_{n \geq 1} &\mapsto (u_1, \cdots, u_n).
 \end{split}
 \end{equation*}
We define the measure $\Pi$ on $c_0(S)$ as follows:
\begin{equation*}
\Pi_{|A_0}(du)= \delta_{A_0}(du),
\end{equation*}
and, for each $n \geq 1$,
\begin{equation*}
\Pi_{|A_n}(du) = \lambda^n (\pi_{n}(du)),
\end{equation*}
where $\lambda^n$ is the Lebesgue measure on $S^n$. We have the following result:
\begin{prpstn}	\label{pamc:prop:abs_cont}
  The law of a Point Absorbed Markov Chain $(M_n)_{n \geq 0}$ on the sphere S, conditioned to start from $a$ is absolutely continuous versus $\Pi$.
\end{prpstn}
\begin{proof}
	Let $\gamma$ be the distribution of $(M_n)_{n \geq 0}$ conditioned to start from $a$.
	We fix $A \in \mathcal{B}(c_0(S))$ such that $\Pi(A)=0$.
	Since $A_0$ is an atom for $\Pi$, $\Pi(A)=0$ implies that $A_0 \cap A = \emptyset$. 
	Thus, we just have to check that $\gamma(A \cap A_n)=0$, for every $n \geq 1$, and to apply the fact that
	\begin{equation*}
 		\gamma(A)= \sum_{n=0}^{+\infty} \gamma(A \cap A_n)=0.
 	\end{equation*}
	Subsequently, using the Markov property, we write:
 	\begin{equation*}
 		\begin{split}
 			\gamma(A \cap A_n)&= \mathbb{P}_a((M_1, \cdots, M_n) \in \pi_{n}(A \cap A_n), M_{n+1} = 0) \\
 			&= \mathbb{P}_a((M_1, \cdots, M_n) \in \pi_{n}(A \cap A_n)) \, P_{M_n}(0) \\
 			&\leq \int_{\pi_{n}(A \cap A_n)} p(a,u_1) \cdots p(u_{n-1},u_{n}) \, du_1 \cdots du_n \\
 		\end{split}
 	\end{equation*}
 	since $P_x(0) \leq 1$, for every $x \in S_0$.
	The desired result follows when we recall that $p$ is absolutely continuous with respect to $\lambda$ and $\lambda(\pi_{n-1}(A \cap A_n))=0$.
\end{proof}
	
This last result allows us to use the Radon-Nykodym-Lebesgue theorem that provide the existence of a density with respect to $\Pi$ for the distribution $\gamma$.
The point is now to exhibit this density.
\begin{prpstn}	\label{pamc:prop:density}
	The density with respect to $\Pi$ of the law of the PAMC on the sphere $(M_n)_{n \geq 0}$ conditioned to start from the point $a\not=0$, is
 	\begin{equation*}
  		P_{a}(0) 1_{A_0}(u) + \sum_{n=1}^{+\infty} p(a,u_1) \cdots p(u_{n-1},u_n) \, 1_{S-\{0\}}(u_1) \cdots 1_{S-\{0\}}(u_n) \, P_{u_n}(0) \, 1_{A_n}(u).
  	\end{equation*}
  	In addition, this density is normalized.
\end{prpstn}
\begin{proof}
	In order to prove this result, we must show, for each Borelian cylinders of finite dimension $C \in \mathcal{F}$, that
		\begin{equation*}
 			\gamma(C)= \int_C P_{a}(0) 1_{A_0}(u) + \sum_{n=1}^{+\infty} p(a,u_1) \cdots p(u_{n-1},u_n) \, 1_{S-\{0\}}(u_1) \cdots 1_{S-\{0\}}(u_n) \, P_{u_n}(0) \, d\Pi(du).
 		\end{equation*}
 	We start by recalling that a Borelian cylinders of finite dimension has the form $C_0 \times \cdots \times C_m \times S \times \cdots$ and the fact that 
 		\begin{equation*}
 		\gamma(C)= \sum_{n=0}^{+\infty} \gamma(C \cap A_n).
 		\end{equation*}
 	since the sequence $(A_n)_{n \geq 0}$ forms a partition of $C_0(S)$.
 	If $n >0$, we can observe that, for $n < m$,
 	 	\begin{equation*}
 	 		\begin{split}
 				\gamma(C \cap A_n)&=\mathbb{P}_a(M_1 \in C_1-\{0\}, \cdots, M_n \in C_n-\{0\}, M_{n+1}=0) \\
 				&= \int_{C_1-\{0\}} \cdots \int_{C_n-\{0\}} p(a,u_1) \cdots p(u_{n-1},u_n) \, P_{u_n}(0) \, du_1 \cdots du_n  \\
 				&= \int_{\pi_n(C \cap A_n)} p(a,u_1) \cdots p(u_{n-1},u_n) \, P_{u_n}(0) \, du_1 \cdots du_n  \\
 				&= \int_{C} p(a,u_1) \cdots p(u_{n-1},u_n) \, P_{u_n}(0) \, 1_{A_n}(u) \, d\Pi(du)  \\
			\end{split} 		
 		\end{equation*}
 	 or, for $n > m$,
 	\begin{equation*}
 	 	\begin{split}
 			\gamma(C \cap A_n)&=\mathbb{P}_a(M_1 \in C_1-\{0\}, \cdots, M_m \in C_m-\{0\}, \cdots ,M_n \in S-\{0\}, M_{n+1}=0) \\
 			&= \int_{C_1-\{0\}} \cdots \int_{C_m-\{0\}} \, \int_{S-\{0\}} \cdots \int_{S-\{0\}} p(a,u_1) \cdots p(x_{n-1},u_n) \, P_{u_n}(0) \, du_1 \cdots du_n  \\
 			&= \int_{\pi_n(C \cap A_n)} p(a,u_1) \cdots p(u_{n-1},u_n) \, P_{u_n}(0) \, du_1 \cdots du_n  \\
 			&= \int_{C} p(a,u_1) \cdots p(u_{n-1},u_n) \, P_{u_n}(0) \, 1_{A_n}(u) \, d\Pi(du)  \\
		\end{split} 		
 	\end{equation*}
 	which show the first part of the result, since 
 	\begin{equation*}
 			\gamma(C \cap A_0) = \mathbb{P}_a(M_1=0) = P_a(0) = \int_{C} P_a(0) \delta_{A_0} = \int_{C} P_a(0) 1_{A_0}(u) d\Pi(du).		
 	\end{equation*}	
As it is not obvious in the proof, we show that the density is normalized using the fact that
 	\begin{equation*}
 		\begin{split}
 			\mathbb{P}(T=n+1) &= \mathbb{P}_a(M_1 \in S-\{0\}, \cdots, M_n \in S-\{0\}, M_{n+1}=0) \\
 			&= \int_{S^{\mathbb{N}}} p(a,u_1) \cdots p(u_{n-1},u_n) P_{u_n}(0) 1_{A_n}(u) \, d\Pi(du), \\
		\end{split}
 	\end{equation*}
and 
 \begin{equation*}
 	\mathbb{P}(T=1) = \mathbb{P}_a(M_1 =0) = \int_{C} P_a(0) 1_{A_0}(u) d\Pi(du), \\
 \end{equation*}
where $T$ is the first time $M$ reaches the absorbing point $0$.
In fact, this is enough when we know that $T$ is almost surely finished and that $\mathbb{P}(T=0)=0$.
\end{proof}

\subsection{A Class of $\Pi$-irreducible Instrumental Kernels} \label{subsection:some:Pi:irreductible}

The Hastings-Metropolis algorithm was originally designed for real random variables and has been widely used in this case.
As a result, extensive studies have been been made to compare different instrumental kernels and show that they play a major role on the reliability of the samples.
Since it is quite new to extend the algorithm to the PAMC on the sphere $S$, we will just give an admissible class of kernels without debates and deep examinations.

The main property required by the Hastings-Metropolis algorithm on an instrumental kernel is the $\gamma$-irreducibility in the sense defined below, as it is a necessary condition for the convergence of the algorithm \cite{meyn, tier}.
Subsequently, we introduce the following probability transition kernel on $c_0(S) \times \mathcal{B}(c_0(S))$ in term of its density with respect to $\Pi$ :
\begin{equation*}
 \kappa(u,dv)= \Theta_0(u) 1_{A_0}(v) + \sum_{k=1}^{+\infty} \Theta_k(u)  \nu_k(u,v) 1_{A_k}(v) 
\end{equation*}
where we assume that:
\begin{enumerate}
 \item For each $u \in c_0$, the sum of the $(\Theta_k(u))_{k \geq 0}$ is 1.
 \item For each $u \in c_0$ and $k \geq 0$, $\Theta_k(u) > 0$.
 \item For each $k \geq 1$, $\nu_k(u,dv)$ is a probability transition kernel on $S_0^k$ having a density with respect to $\lambda^k$.
\end{enumerate}
This statement ensures that $\kappa$ is a probability transition kernel on $c_0$. 
We describe the behavior of the chain:
\begin{enumerate}
 \item We change the number of non-null points using the family $(\Theta_k(u))_{k \geq 0}$.
 	For example, suppose that $u \in A_m$, then $u$ moves into $A_k$ with the probability $\Theta_k(u)$.
 	As a result, the chain $u$ loses or gains points different from $0$.
 	In the case consisting of adding new points, we choose to initialize all of them at a position $b \in S-\{0\}$.
 	Otherwise, by losing, we mean that the last $m-k$ positions are set to 0.
 \item We use a classical instrumental kernel on the finite dimensional vector of non-null positions.
\end{enumerate} 
Before proving any property on this kernel, we give a set of definitions to understand the concept of the irreducibility of a Markov chain: 
\begin{dfntn}
 	Let $G$ be a topological space, $\mathcal{G}$ a $\sigma$-algebra on $G$, $m$ a probability measure and $\mu$ a probability transition kernel.
 	We say that $A \in \mathcal{G}$ is attainable from $x \in G$ if:
 	\begin{equation*}
  		\text{it exists } n > 1 \text{ such that }\mu^n(u,A) > 0,
 	\end{equation*}
 	and attainable from $x \in G$ in one step if $\mu(u,A) > 0$.
\end{dfntn}
\begin{dfntn}
	Let $G$ be a topological space, $\mathcal{G}$ a $\sigma$-algebra on $G$, $m$ a probability measure and $\mu$ a probability transition kernel.
		\begin{enumerate}
		\item $B \in \mathcal{G}$ is $m$-communicating if
 			\begin{equation*}
  				\forall x \in B, \forall A \in \mathcal{G}\text{ such that }A \subset B \text{, } m(A)>0, A\text{ is attainable from }x.
			\end{equation*}
		\item $B \in \mathcal{G}$ is quickly $m$-communicating if
			\begin{equation*}
  				\forall x \in B, \forall A \in \mathcal{G}\text{ such that }A \subset B \text{, } m(A)>0, A\text{ is attainable in one step from }x.
			\end{equation*}
	\end{enumerate}
\end{dfntn}
\begin{dfntn}
	Let $G$ be a topological space, $\mathcal{G}$ a $\sigma$-algebra on $G$, $m$ a probability measure and $\mu$ the probability transition kernel of a Markov chain $(X_n)_{n \geq 0}$.
	\begin{enumerate}
		\item $G$ is $m$-communicating, $(X_n)_{n \geq 0}$ and $\mu$ are said $m$-irreducible.
		\item $G$ is quickly $m$-communicating, $(X_n)_{n \geq 0}$ and $\mu$ are said strongly $m$-irreducible.
	\end{enumerate}
\end{dfntn}
From \cite{tier}, we know that: if $\kappa$ is $\Pi$-irreducible, then $\kappa$ is also $\gamma$-irreducible since $\gamma$ is absolutely continuous with respect to the measure $\Pi$.
As a result, the result that follows provides the property required for the convergence which is aforementioned.
\begin{prpstn} \label{instru:prop:irreduc}
If $\kappa$ is such that, for each $k \geq 1$, $\nu_k(u,dv)$ is strongly $\lambda^k$-irreducible.
Then, $\kappa$ is strongly $\Pi$-irreducible.
\end{prpstn}
\begin{proof}
Let $A \in \mathcal{B}(c_0)$ be a $\Pi$-positive subset and $u \in c_0$ a sequence. 
In order to prove that $\kappa$ is strongly $\Pi$-irreducibility, we have to show that $\kappa(u,A) > 0$.
Note that this result holds if, for each $k \geq 0$,
\begin{equation*}
	 A \subset A_k \text{ }\text{ is } \Pi \text{-positive} \implies \text{ A is attainable from } u \in c_0.
\end{equation*}
Let fix $k \geq 0$ and assume that $A \subset A_k$. 
From the definition of $\kappa$, we have
 \begin{equation*}
  \kappa(u,A) = \int_A \Theta_k(u) \nu_{k}(u,dv) \,dv.
 \end{equation*} 
Since $\Theta_k(u)>0$ for every $ k > 0$ and $u \in c_0(S)$, we only have to prove that 
 \begin{equation*}
  \int_A  \nu_{k}(u,v) \,dv > 0.
 \end{equation*}
 The absolute continuity and the fact that $\nu_k(u,v)$ is strongly $\lambda^k$-irreducible induce that
 \begin{equation*}
    \text{ if } A \text{ is a }\lambda^k \text{-positive set}, \text{ then } \nu_k(u,A) > 0.
 \end{equation*}
 and the result holds.
 Indeed, if we suppose the opposite, then we have a conflict with the strong $\lambda^k$-irreducibility, since
  \begin{equation*}
    \text{ if } A \subset A_k \text{ is a }\Pi \text{-positive set}, \text{ then } A \text{ is a }\lambda^k \text{-positive set}.
 \end{equation*}
\end{proof}

\subsection{Convergence of the extended Hastings-Metropolis Algorithm} \label{subsection:convergence:HM}

Before the proof of convergence of the algorithm, we give an example of $(\Theta_k(u))_{k\geq 0}$ and $(\nu_k(u,v))_{k\geq 1}$ such that $\kappa$ is $\gamma$-irreducible.
Let $G$ be the shifted geometric distribution on $\mathbb{N}$ and $g$ the density of the uniform distribution on $S_0$. 
For each $u \in c_0(S)$, we set $\Theta_k(u)=\mathbb{P}(G=k)$ and
\begin{equation*}
 \nu_k(u,v)=\prod_{i=1}^k g(v).
\end{equation*}

The following theorem is the main theoretical result of this paper. 
It relies on the topological claim which provides the hypothesis required in the theoretical results used for the proof.
We decide to present a theorem with relatively strong hypothesis in order to convince the reader of the convergence of the more complex case used in the numerical experiments.
\begin{thrm}
	Let $(M_n)_{n \geq 0}$ be a Point Absorbed Markov chain on $S$ starting from $a$.
	We consider the following instrumental kernel:
	\begin{equation*}
 		\kappa(u,dv)= \Theta_0(u) 1_{A_0}(v) + \sum_{k=1}^{+\infty} \Theta_k(u)  \nu_k(u,v) 1_{A_k}(v) 
	\end{equation*}
	satisfying the following hypothesis:
	\begin{enumerate}
 		\item For each $u \in c_0$ and $k \geq 0$, 
 			\begin{equation*}
 				\Theta_k(u) = \mathbb{P}(G=k)
 			\end{equation*}
 			where:
 			\begin{enumerate}
 				\item $G$ a probability law on $\mathbb{N}$.
 				\item for every $k \in \mathbb{N}$, $\mathbb{P}(G=k) > 0$.
 			\end{enumerate}
 		\item for each $k \geq 1$, 
 		\begin{equation*}
 			\nu_k (u,v) = h(u_1,v_1) \cdots h(u_k,v_k)
 		\end{equation*}
 		where:
 		\begin{enumerate}
 			\item $h$ is a probability transition kernel on $S_0$.
 			\item $h$ is absolutely continuous with respect to $\lambda$.
 			\item $h$ is strongly $\lambda$-irreducible.
 			\item $h$ is symmetric: $h(x,y)=h(x,y)$, for every $(u,v)$ in $S_0^2$.
 			\end{enumerate}
	\end{enumerate}
Then, the Hastings-Metropolis kernel $\Gamma$ formed with $\kappa$ and $r$ converges to $\gamma$ with respect to the topology of the total variation norm.
\end{thrm}
\begin{proof}
The probability transition kernel $h$ on $S$ is strongly $g$-irreducible, since it is absolutely continuous with respect to $\lambda$ and strongly $\lambda$-irreducible.
In addition, $\nu_k (u,v)$ is strongly $\lambda^k$-irreducible as a product of strongly $\lambda$-irreducible kernel . 
Using Proposition \ref{instru:prop:irreduc}, we conclude that the kernel $\kappa$ is strongly $\gamma$-irreducible.

In order to prove the convergence the Hastings-Metropolis kernel $\Gamma$, we follow \cite{tier} which shows that we just have to show that $\Gamma$ is $\gamma$-irreducible and $\gamma\{\overline{r}(u)>0\}>0$ to obtain the convergence with respect to the topology of the total variation norm.
Before starting the proof, we recall that 
\begin{equation*}
\Gamma(u,dv)=\overline{\kappa}(u,v)\Pi(dv)+\overline{r}(u)\delta_u(dv) 
\end{equation*}
where
\begin{equation*}
 \overline{\kappa}(u,v)=
 \begin{cases}
  \kappa(u,v) r(u,v),\quad &\text{if } u\not= v, \\
  0,&\text{if } u=v,
 \end{cases}
\end{equation*}
\begin{equation*}
 \overline{r}(u)=1-\int \overline{\kappa}(u,v) \Pi(dv).
\end{equation*}
and 
\begin{equation*}
 r(u,v)=
 \begin{cases}
  \min\left\lbrace \dfrac{\gamma(v) \kappa(v,u)}{\gamma(u) \kappa(u,v)},1 \right\rbrace, &\text{if }\gamma(u) \kappa(u,v)>0\\
  1, &\text{if }\gamma(u) \kappa(u,v)=0
 \end{cases}.
\end{equation*}

We start by showing that $\Gamma$ is strongly $\gamma$-irreducible. Let $A \in \mathcal{B}(c_0)$ be a $\gamma$-positive subset and $u$ a sequence of $c_0(S)$ such that $u \in A_l$. 
We can establish that $\Gamma$ is strongly $\gamma$-irreductible if we prove that $\Gamma^n(u,A) > 0$.
We use the approach developed in the proof of Proposition \ref{instru:prop:irreduc}.
Let fix $k \geq 0$ and suppose that $A \subset A_k$. 
With the second term in the expression of $\Gamma$ and the fact that
\begin{equation*}
	r(u,v) \geq \dfrac{\gamma(v) \kappa(v,u)}{\gamma(u) \kappa(u,v)},
\end{equation*}
it is enough to show that
\begin{equation*}
	\begin{split}
 		\overline{\kappa}(u,A) &= \int_A \kappa(u,v) \dfrac{\gamma(v) \kappa(v,u)}{\gamma(u) \kappa(u,v)} \Pi(dv) \\
 		 &= \int_A \dfrac{\gamma(v) \kappa(v,u)}{\gamma(u)} \lambda^k(dv)> 0.
	\end{split}
\end{equation*}
Moreover, we can suppose that $u$ is such that $\gamma(u) > 0$ on $B \subset A$, else the result is proved since
\begin{equation*}
	\begin{split}
 		\overline{\kappa}(u,A) &= \int_A \kappa(u,v) \Pi(dv) \\
	\end{split}
\end{equation*}
and $\kappa$ is strongly $\gamma$-irreducible.
As result, we have  
\begin{equation*}
	\begin{split}
		\int_A \dfrac{\gamma(v) \kappa(v,u)}{\gamma(u)} \lambda^k(dv) \geq \frac{1}{\gamma(u)} \, \int_B \gamma(v) \kappa(v,u) \lambda^k(dv).
	\end{split}
\end{equation*}
Thereupon, we can suppose that it exists $B' \subset B$ such that $\gamma(v) > C_1$, for each $v \in B'$, since $\gamma(A) > 0$. 
Thus, we get that 
\begin{equation*}
	\begin{split}
		\int_A \dfrac{\gamma(v) \kappa(v,u)}{\gamma(u)} \lambda^k(dv) \geq \frac{C_1}{\gamma(u)} \, \int_{B'} \Theta_l(v) \, h(v_1,u_1) \cdots h(v_n,u_n) \, dv_1 \cdots dv_n.
	\end{split}
\end{equation*}
Using the symmetry of $h$, we can rewrite:
\begin{equation*}
	\begin{split}
		\int_A \dfrac{\gamma(v) \kappa(v,u)}{\gamma(u)} \lambda^k(dv) \geq \frac{C_1}{\gamma(u)} \, \int_{B'} \Theta_l(v) \, h(u_1,v_1) \cdots h(u_n,v_n) \, dv_1 \cdots dv_n \geq \frac{C}{\gamma(u)} \, \inf_{v \in c_0(S)}{\Theta_l(v)},
	\end{split}
\end{equation*}
and the strong $\gamma$-irreducibiblity follows from the hypothesis. 

The last step consist of showing that $\gamma\{\overline{r}(u) > 0\}>0$.
From the definition of a PAMC, for every $l \geq 0$, we know that $\gamma(A_l)>0$.
Suppose $v \in A_k$ with $k < l$. Then, for every $u \in A_l$, we have:
\begin{equation*}
	\frac{\gamma(v) \kappa(v,u)}{\gamma(u) \kappa(u,v)}=\frac{p(a,v_1) \cdots p(v_{k-1},v_k) P_{v_k}(0)}{p(a,u_1) \cdots p(u_{l-1},u_l) P_{u_l}(0)} \times \frac{\Theta_l(v)}{\Theta_k(u)}
	\times \frac{h(v_1,u_1) \cdots h(v_l,u_l)}{h(u_1,v_1) \cdots h(u_k,v_k)}.
\end{equation*}
Since $h$ is symmetric, we can rewrite:
\begin{equation*}
	\frac{\gamma(v) \kappa(v,u)}{\gamma(u) \kappa(u,v)}=\frac{p(a,v_1) \cdots p(v_{k-1},v_k) P_{v_k}(0)}{p(a,u_1) \cdots p(u_{l-1},u_l) P_{u_l}(0)} \times \frac{\Theta_l(v)}{\Theta_k(u)}
	\times h(v_{k+1},u_{k+1}) \cdots h(v_l,u_l).
\end{equation*}
Moreover, $h(v,\cdot)$ being absolutely continuous with respect to the Lebesgue measure on $S$, $h(v,.)$ is continuous and bounded on $S$.
But, $h$ is symmetric. Thus, $h$ is uniformly bounded on $S$ and 
\begin{equation*}
	\frac{\gamma(v) \kappa(v,u)}{\gamma(u) \kappa(u,v)} \leq C_2 \, \frac{p(a,v_1) \cdots p(v_{k-1},v_k) P_{v_k}(0)}{p(a,u_1) \cdots p(u_{l-1},u_l) P_{u_l}(0)} \times \frac{\Theta_l(v)}{\Theta_k(u)}.
\end{equation*}
From the hypothesis on the transition kernel of a PAMC, we have:
\begin{equation*}
	\frac{\gamma(v) \kappa(v,u)}{\gamma(u) \kappa(u,v)} \leq \frac{C_2}{c} \times \frac{M^k}{m^l} \times \frac{\Theta_l(v)}{\Theta_k(u)}.
\end{equation*}
Since the sum of the $\Theta_l(v)$ is finite for every $v \in c_0(S)$ and the same, the sequence $(\Theta_l(v))_{l \geq 0}$ converges to $0$.
As a result, we can choose $l$ such that 
\begin{equation*}
	\Theta_l(v) < \left( \frac{C_2}{c} \times \frac{M^k}{m^l} \times \frac{1}{\Theta_k(u)} \right)^{-1}.
\end{equation*}
Thus, for every $u \in A_l$ and $v \in A_k$,
\begin{equation*}
	\frac{\gamma(v) \kappa(v,u)}{\gamma(u) \kappa(u,v)} < 1,
\end{equation*}
and, for every $u \in A_l$,
\begin{equation*}
	\int \overline{\kappa}(u,v) \, \Pi(dv) \leq  \int_{c_0(s)-\{A_k\}} \kappa(u,v) \,\Pi(dv) 
	+ \int_{A_k} \kappa(u,v) \, \frac{\gamma(v) \kappa(v,u)}{\gamma(u) \kappa(u,v)}\Pi(dv) < 1
\end{equation*}
which provides the desired results.
\end{proof}

\section{Practical implementation for the monokinetic particle simulation} \label{section:one:two:dim:case}

We present a model of a monokinetic particle that travels along straight lines between random collision points. The  sequence of collision points constitutes a Markov chain which is almost surely absorbed in finite time. This Markov chain is identical to the Boundary Absorbed Markov Chain of Section \ref{section:HM:on:Markov:chains}, except that we consider here that absorption can occur with positive probability not only outside of a domain but at any point of space. I short words, in Section \ref{section:HM:on:Markov:chains}, the Markov chain is killed when it leaves the domain while in Section \ref{section:one:two:dim:case}, it is killed when it leaves the domain and also possibly in the domain.
Thus, the notion of pdf for the space of monokinetic particle trajectories must be first defined, in a different way than in Section \ref{section:HM:on:Markov:chains}. Then, we present one and two-dimensional versions of the monokinetic particle model, the instrumental kernels we consider, and we give the corresponding explicit expressions of the unconditional and conditional pdf of the trajectories. The final version of Algorithm \ref{algo:AG:with:HM} for the monokinetic particle simulation is then summed up.

\subsection{General vocabulary and notation} \label{subsection:physical:vocabulary}

Throughout Section \ref{section:one:two:dim:case}, we consider a monokinetic particle (a particle with constant speed and yielding no subparticle birth) evolving in $\xR^d$, with $d=1,2$.
The birth of the particle takes place at $s$, which we write as $X_0 = s$.
Then, the trajectory of the monokinetic particle is characterized by its collision points, which constitute a homogeneous Markov chain $(X_n)_{n \in \xN^*}$ on $\xR^d \cup \{\Delta\}$ with transition kernel
\begin{equation} \label{eq:homogeneous:transition:monokinetic}
k(x_{n},dx_{n+1}) = \delta_\Delta(dx_{n+1}) \mathbf{1} {\{x_{n} = \Delta\}}  +  \left\{ P(x_n) \delta_\Delta(dx_{n+1}) + \left[1 - P(x_n)\right] q(x_n,x_{n+1}) dx_{n+1}  \right\} \mathbf{1} {\{x_{n} \neq \Delta\}}.
\end{equation}
In the above display,
$P(x_n) \in [0,1]$ is the probability of absorption for a collision taking place at $x_n \in \xR^d$. 
Absorption at collision $n$ is here conventionally defined as $X_{n+1} = \Delta = 0 \in \xR^{d+1}$, which implies
$X_m = \Delta$ for any $m > n$. We call $\Delta$ the absorbed state, or resting point as in Section \ref{section:HM:on:Markov:chains} and use the convenient convention that an absorbed monokinetic particle makes an infinite number of collisions at $\Delta$. Finally, conditionally to the collision $x_n$, the particle is scattered with probability $1-P(x_n)$, in which case the next collision point has pdf $q(x_n,.)$.

We assume here, similarly to Section \ref{section:HM:on:Markov:chains}, that the Markov chain $(X_n)_{n \in \xN^*}$ has the property that absorption happens almost surely after a finite number of collisions. That is, almost surely, there exists $m \in \xN$ so that $X_n = \Delta$ for $n \geq m$. 
This assumption holds for example when $P(x_n) = 1$ out of a compact set $C$ of $\xR^d$ and where there exists a positive constant $c$ so that $q(x,\mathbb{R}^d \backslash C) \geq c$ for all $x \in \mathbb{R}^d$, which is the case in Section \ref{section:one:two:dim:case}.
We say that the monokinetic particle is active at time $n$, or at $X_n$, or before collision $n$, if $X_n \neq \Delta$.

Finally, note that the Markov Chain of the collision points $(X_n)_{n \in \xN^*}$ does not include the birth point $X_0 = s$, which entails no loss of information since $s$ is deterministic.

\subsection{The measured space of monokinetic particle trajectories} \label{subsection:base:measure:monokinetic}

For further reference throughout Section \ref{section:one:two:dim:case}, we define here the measured space $(c_0,\mathcal{S},\Pi)$ of the monokinetic particle paths.
We start by defining $c_0$ and $\mathcal{S}$.
\begin{dfntn} \label{def:tribu:markov:chain}
Define
\[
c_0=\{(x_n)_{n \geq 1} \in \left( \xR^d \cup \{\Delta\} \right)^{\mathbb{N}^*}: \exists n_0 \in \mathbb{N}^*, \forall n < n_0, x_n \neq \Delta; \forall n \geq n_0, x_n= \Delta \}.
\]
Let $\mathcal{S}$ be the smallest sigma-algebra on $c_0$ containing the sets
$\{x \in c_0 | x_{1} \in B_{1}  , ... ,x_{n} \in B_{n}, x_{n+1} = \Delta \}$, for $n \in \xN$ and
$B_i \in \mathcal{B} \left(\xR^d\right)$, where $\mathcal{B} \left( \xR^d \right)$ is the Borel sigma-algebra on $\xR^d$.
\end{dfntn}

We define for $n \geq 0$
\begin{equation} \label{eq:An:in:Section:four:five}
A_n = \{ x \in c_0 ; \forall \, 1 \leq j \leq n: x_j \neq \Delta, \forall k \geq n+1: x_k = \Delta  \},
\end{equation}
that is the set of trajectories that are absorbed at collision point $n$ (so that they are in the absorbed state from collision point $n+1$ and onward). Note that the $A_n$, for $n \geq 0$, constitute a partition of $c_0$.
The existence of the measure $\Pi$ is now shown in the following proposition, which can be proved in the same way as in Section \ref{section:HM:on:Markov:chains}.
\begin{prpstn} \label{prop:measure_on_S}
There exists a unique measure $\Pi$ on
$(c_0 , \mathcal{S})$ that verifies the following relation, for any $E_n = \{ x \in A_{n};  x_1 \in B_1 ,...,x_n \in B_n \}$, with $B_1,...,B_n \in \mathcal{B} \left(\xR^{d} \right)$ and $n \in \xN$:
\begin{equation} \label{eq: def_measure_A}
\Pi(E_n) =
\lambda (B_1)... \lambda(B_n),
\end{equation}
with $\lambda$ the Lebesgue measure on $\xR^d$, and with $\Pi(E_0) = \Pi(A_0) =0$.
\end{prpstn}

\subsection{Description of the one-dimensional case and expression of the probability density functions} \label{subsection:one:dim:case}


\subsubsection{A one-dimensional random walk} \label{subsection:one:dim:procedure}

We consider that the monokinetic particle evolves in $\xR$. With the notation of \eqref{eq:homogeneous:transition:monokinetic}, we set $q(x_n,.)$ as the Gaussian pdf with mean $0$ and variance $\sigma^2$ and we set $P(t) = \mathbf{1} {\{t \not \in (A,B)\}} + P \mathbf{1} {\{t \in (A,B)\}}$, with $A < 0 < B $ and $0 < P < 1$. Thus, the particle travels with normally distributed increments, has a probability of absorption $P$ at each collision point in the domain of interest $(A,B)$ and is absorbed if it leaves this domain.

The following algorithm, when tuned with source point $s=0$, sums up how one can sample one-dimensional trajectories. 

\begin{lgrthm} \label{algo:sampling:f:one:dim}

Objective: from a source point $s \in \mathbb{R}$ and the parameters $\mathcal{D} = (A,B)$, $\sigma^2$ and $P$, sample a trajectory $x$ as described above.

~

\begin{itemize}
\item Set $i=0$, $x_i=s$ and``state = active''.
\item While ``state = active'' do
\begin{itemize}
\item Sample $x_{i+1}$ from the $\mathcal{N}(x_i,\sigma^2)$ distribution.
\item If $x_{i+1 } \not \in \mathcal{D}$
\begin{itemize}
\item Set ``state = inactive''.
\end{itemize}
\item If $x_{i+1 } \in \mathcal{D}$
\begin{itemize}
\item With probability $P$, set ``state = inactive''.
\end{itemize}
\item Set $i=i+1$.
\end{itemize}
\item Return the infinite sequence $(x_1,...,x_{i},\Delta,...)$.
\end{itemize}

\end{lgrthm}

The event of interest is here that the monokinetic particle reaches the domain $(- \infty,A ]$. When using the interacting-particle method of Section \ref{section:reminder:AG}, this event is expressed by $\Phi(x) \geq 0$, with $\Phi(x) = A - \inf_{ i \in \xN^*; x_i \neq \Delta} x_i$. Note that, almost-surely, the infimum is taken over a finite number of points.

Although the two-dimensional case of Section \ref{subsection:two:dim:case} is more realistic, we address here
absorption with positive probability at each collision point, which is an important features of shielding studies by Monte Carlo code. 
Furthermore, by setting $P$ sufficiently large, and $A$ sufficiently away from $0$, we will see that we can tackle problems of estimation of arbitrary small probabilities. In Section \ref{section:numerical:results}, we will consider a probability small enough so that the interacting-particle method of Section \ref{section:reminder:AG} outperforms a simple Monte Carlo method.

\subsubsection{Expression of the probability density function of a trajectory}

We now give the expression of the pdf (with respect to the setting of Definition \ref{def:tribu:markov:chain} and Proposition \ref{prop:measure_on_S}) of a trajectory obtained from the one-dimensional model above.
We let $(x_i)_{i \in \xN^*}$ be the sequence of collision points (the trajectory) of a monokinetic particle. We let $\mathcal{D} = (A,B)$. We denote $\phi(m,\sigma^2,t)$ the pdf at $t$ of the one-dimensional Gaussian distribution with mean $m$ and variance $\sigma^2$.

\begin{prpstn} \label{prop:pdf:initial:one:dimensional}
The pdf, with respect to $(c_0,\mathcal{S},\Pi)$ of Definition \ref{def:tribu:markov:chain} and Proposition \ref{prop:measure_on_S}, of a trajectory $(x_n)_{n \in \xN^*}$, sampled from Algorithm \ref{algo:sampling:f:one:dim}, is $f(x) = \sum_{n \in \xN^*} \mathbf{1}_{A_{n}}(x) f_n(x)$, with
\[
f_n(x) = \left( \prod_{i=1}^{n-1}    \phi(x_{i-1},\sigma^2,x_i) (1-P) \mathbf{1} {\{x_i \in \mathcal{D}\}} \right) \phi(x_{n-1},\sigma^2,x_n) \left( \mathbf{1} {\{x_n \not \in \mathcal{D}\}} + P \mathbf{1} {\{x_n \in \mathcal{D}\}} \right),
\]
where $x_0=0$ by convention.
\end{prpstn}

The pdf of Proposition \ref{prop:pdf:initial:one:dimensional} has to be evaluated for each trajectory, either sampled from its initial distribution, or from an instrumental kernel $\kappa$ in Algorithm \ref{algo:HM}. The perturbation methods, defining $\kappa$, are presented below.

In Proposition \ref{prop:pdf:initial:one:dimensional}, note that, in case the monokinetic particle leaves the domain $\mathcal{D}$, we explicitly use the exact position of the collision point outside $\mathcal{D}$.
This exact value is not needed to assess if the monokinetic has reached the domain $( - \infty , A]$.
Thus, we might add some variance in the HM method, because we use a source of randomness (the exact collision point at which the monokinetic particle leaves $D$) that does not impact the event of interest. This is nevertheless inevitable if one requires explicit evaluation of the pdf. Indeed, modifying the definition of trajectories and of pdf so that collision points out of the domain are not stored would add to the pdf expression the probability that, starting from a birth or scattering point in the domain $D$, the next collision point lies outside $D$. This probability has an explicit expression in this one-dimensional case, but not in the framework of Section \ref{subsection:two:dim:case}, and {\it a fortiori}
not in shielding studies involving more complex Monte Carlo codes. Thus, to avoid evaluating this probability numerically each time a pdf of a trajectory is computed, we store the collision points outside the domain $D$.

The evaluation of a pdf like that of Proposition \ref{prop:pdf:initial:one:dimensional} is an intrusive operation on a Monte Carlo code. Indeed, it necessitates to know all the random-quantity sampling that are done when this code samples a monokinetic-particle trajectory. Thus, the Monte Carlo code is not used as a black box.
Nevertheless, the computational cost of the pdf evaluation is of the same order as the computational cost of a trajectory sampling, and the same kind of operations are involved. Namely, both tasks require a loop which length is the number of collisions made by the monokinetic-particle before its absorption. Furthermore, for each random quantity that is sampled for a trajectory sampling, the pdf evaluation requires to compute the corresponding pdf. For example, in the case of
Proposition \ref{prop:pdf:initial:one:dimensional}, when a trajectory sampling requires to sample $n$ Gaussian variables and $n$ or $n-1$ Bernoulli variables, the trajectory-pdf evaluation requires to compute the corresponding Gaussian pdf and Bernoulli probabilities.

Finally, the discussion above holds similarly for the two-dimensional case of Section \ref{subsection:two:dim:case}.

\subsubsection{Description of the trajectory perturbation method when $P = 0$} \label{subsubsection:one:dim:perturbation:procedure}

For clarity of exposition, we present first the perturbation method when $P = 0$.
In this case, the monokinetic particle is a random walk on $\xR$, that is absorbed once it goes outside $\mathcal{D}$.

The perturbation method is parameterized by $\tilde{\sigma}^2 > 0$.
Let us consider a historical trajectory $(x_i)_{i \in \xN^*}$, absorbed at collision $n$. Then, the set of birth and collision points of the perturbed monokinetic-particle is an inhomogeneous Markov chain $(Y_i)_{i \in \xN}$ so that $Y_0 = 0$. If $i \leq n-1$, and if the perturbed monokinetic particle is still in $\mathcal{D}$ at collision point $i$, we have $Y_{i+1} = Y_i + \epsilon_{i+1}$, where the
$(\epsilon_i)_{1 \leq i \leq n}$ are independent and where $\epsilon_i$ follows a $\mathcal{N}(x_i - x_{i-1},\tilde{\sigma}^2)$ distribution.

Similarly to the initial sampling, the perturbed monokinetic particle is absorbed at the first collision point outside $\mathcal{D}$.
If the collision point $Y_n$ of the perturbed monokinetic particle is in $\mathcal{D}$ (contrary to $x_n$ for the initial trajectory), the sequel of the trajectory of the perturbed monokinetic particle is sampled as the initial monokinetic particle would be sampled if its collision point $n$
was $Y_n$.

This conditional sampling method for perturbed trajectories is intrusive: it necessitates to change the stochastic dynamic of the monokinetic particle. Nevertheless, the new dynamic is here chosen as to have the same cost as the unconditional sampling, and to require the same type of computations. This is similar to the discussion following Proposition \ref{prop:pdf:initial:one:dimensional}.

\subsubsection{Expression of the probability density function of a perturbed trajectory when $P=0$} \label{subsubsection:perturbed:pdf:no:abs}

\begin{prpstn} \label{prop:pdf:perturbed:one:dimensional:Pa:zero}
Let us consider a historical trajectory $(x_i)_{i \in \xN^*}$, absorbed at collision $n$.
The conditional pdf, with respect to $(c_0,\mathcal{S},\Pi)$ of Definition \ref{def:tribu:markov:chain} and Proposition \ref{prop:measure_on_S}, of a trajectory $(y_n)_{n \in \xN^*}$ sampled from the procedure of Section \ref{subsubsection:one:dim:perturbation:procedure}, is $\kappa(x,y) = \sum_{m \in \xN^*} \mathbf{1}_{A_{m}}(y) f_{n,m}(x,y)$ where, if $m \leq n $
\[
f_{n,m}(x,y) = \prod_{i=1}^{m-1} \left(   \phi(y_{i-1} + ( x_i - x_{i-1}),\tilde{\sigma}^2,y_i) \mathbf{1} {\{y_i \in \mathcal{D}\}} \right)
 \phi(y_{m-1} + ( x_{m} - x_{m-1} ),\tilde{\sigma}^2,y_{m})  \mathbf{1} {\{y_{m} \not \in \mathcal{D}\}},
\]
and if $m > n $,
\begin{eqnarray*}
f_{n,m}(x,y) & = &
\prod_{i=1}^{n} \left(   \phi(y_{i-1} + ( x_i - x_{i-1}),\tilde{\sigma}^2,y_i) \mathbf{1} {\{y_i \in \mathcal{D}\}} \right) \\
&  & \times \prod_{i=n+1}^{m-1} \left(   \phi(y_{i-1} ,\sigma^2,y_i) \mathbf{1} {\{y_i \in \mathcal{D}\}} \right) \\
& &
 \phi(y_{m-1} ,\sigma^2,y_{m})  \mathbf{1} {\{y_{m} \not \in \mathcal{D}\}},
\end{eqnarray*}
where $x_0 = y_0 = 0$ by convention.
\end{prpstn}

Similarly to the discussion following \ref{prop:pdf:initial:one:dimensional}, the computation of the conditional pdf of a perturbed trajectory has the same computational cost as the sampling of this perturbed trajectory.

\subsubsection{Description of the trajectory perturbation method when $P > 0$} \label{subsubsection:one:dim:perturbation:procedure:Pa:non:zero}

In the general case where $P > 0$, the perturbation method is parameterized by $\tilde{\sigma}^2 > 0$ and $0 < Q < 1$.
Let us consider a historical trajectory $(x_i)_{i \in \xN^*}$, absorbed at collision $n$.
As when $P = 0$, the set of birth and collision points of the perturbed monokinetic particle is an inhomogeneous Markov chain $(Y_i)_{i \in \xN}$, so that $Y_0 = 0$.
As when $P = 0$, we modify the increments of the initial trajectory, and, if the perturbed trajectory outsurvives the initial one, we generate the sequel with the initial distribution. Specifically to this case $P > 0$, we perturb the absorption/non-absorption sampling by changing the initial values with probability $Q$.

More precisely, for $i \leq n-1$ and if the perturbed monokinetic particle has not been absorbed before collision point $i$, it is absorbed with probability $max(Q,\mathbf{1} {\{Y_i \not \in D\}})$. If it is scattered instead, we have $Y_{i+1} = Y_i + \epsilon_{i+1}$, where the
$(\epsilon_i)_{1 \leq i \leq n}$ are independent and where $\epsilon_i$ follows a $\mathcal{N}(x_i - x_{i-1},\tilde{\sigma}^2)$ distribution.
If the perturbed monokinetic particle has not been absorbed before collision point $n$, then it is absorbed if $Y_n \not \in \mathcal{D}$. If $Y_n  \in \mathcal{D}$, the perturbed monokinetic particle is absorbed with probability $(1 - Q) \mathbf{1} {\{x_n \in \mathcal{D}\}} + P \mathbf{1} {\{x_n \not \in \mathcal{D}\}}$.
As when $P = 0$, if the perturbed monokinetic particle has not been absorbed before collision point $Y_n$, the sequel of the trajectory of the perturbed monokinetic particle is sampled as the initial particle would be sampled if its collision point $n$ was $Y_n$.

The idea is that, by selecting the difference between $Q$ and $\min(P,1-P)$, the closeness between the perturbed and initial trajectories can be specified, from the point of view of the absorption/non-absorption events.
Finally, the following algorithm sums up how perturbed trajectories can be sampled.
\begin{lgrthm} \label{algo:sampling:kappa:one:dim}

Objective: from an initial trajectory $x$ absorbed at collision $n$ and from the parameters $\mathcal{D} = (A,B)$, $\sigma^2$, $P$, $\tilde{\sigma}^2$ and $Q$, sample a perturbed trajectory $y$ as described above.

~

\begin{itemize}
\item Set $i=0$, $y_i=0$ and ``state = active''.
\item While ``state = active'' and $i+1 < n$ do
\begin{itemize}
\item Sample $y_{i+1}$ from the $\mathcal{N}(y_i + x_{i+1} - x_i,\tilde{\sigma}^2)$ distribution.
\item If $y_{i+1 } \not \in \mathcal{D}$
\begin{itemize}
\item Set ``state = inactive''.
\end{itemize}
\item If $y_{i+1 } \in \mathcal{D}$
\begin{itemize}
\item With probability $Q$, set ``state = inactive''.
\end{itemize}
\item Set $i=i+1$.
\end{itemize}
\item If ``state = inactive'', stop the algorithm and return the infinite sequence $(y_1,...,y_{i},\Delta,...)$.
\item If ``state = active'' do
\begin{itemize}
\item Sample $y_{i+1}$ from the $\mathcal{N}(y_i + x_{i+1} - x_i,\tilde{\sigma}^2)$ distribution.
\item If $y_{i+1 } \not \in \mathcal{D}$
\begin{itemize}
\item Set ``state = inactive''.
\end{itemize}
\item If $y_{i+1 } \in \mathcal{D}$
\begin{itemize}
\item With probability $(1-Q) \mathbf{1} {\{x_{i+1} \in \mathcal{D}\}}  + P \mathbf{1} {\{x_{i+1} \not \in \mathcal{D}\}} $, set ``state = inactive''.
\end{itemize}
\item Set $i=i+1$.
\end{itemize}
\item If ``state = inactive'', stop the algorithm and return the infinite sequence $(y_1,...,y_{i},\Delta,...)$.
\item If ``state = active'' do,
\begin{itemize}
\item  Apply Algorithm \ref{algo:sampling:f:one:dim}, with $s = y_i$ and write $(\tilde{x}_1,...,\tilde{x}_q,\Delta,...)$ for the resulting trajectory.
\item Return the infinite sequence $(y_1,...,y_{i},\tilde{x}_1,...,\tilde{x}_q,\Delta,...)$.
\end{itemize}

\end{itemize}

\end{lgrthm}

\subsubsection{Expression of the probability density function of a perturbed trajectory when $P>0$} \label{subsubsection:one:dim:perturbed:pdf:abs}

\begin{prpstn} \label{prop:pdf:perturbed:one:dimensional:Pa:non:zero}
Let us consider a historical trajectory $(x_i)_{i \in \xN^*}$, absorbed at collision $n$.
Let $y_0 = x_0 = 0$ by convention.
The conditional pdf, with respect to $(c_0,\mathcal{S},\Pi)$ of Definition \ref{def:tribu:markov:chain} and Proposition \ref{prop:measure_on_S}, of a trajectory $(y_n)_{n \in \xN^*}$ sampled from Algorithm \ref{algo:sampling:kappa:one:dim}, is $\kappa(x,y) = \sum_{m \in \xN^*} \mathbf{1}_{A_{m}}(y) f_{n,m}(x,y)$ where, if $m \leq n-1 $,
\begin{eqnarray*}
f_{n,m}(x,y) & = & \prod_{i=1}^{m-1} \left(    \phi(y_{i-1} + ( x_i - x_{i-1}),\tilde{\sigma}^2,y_i) (1 - Q) \mathbf{1} {\{y_i \in \mathcal{D}\}} \right) \\
& & \phi(y_{m-1} + ( x_{m} - x_{m-1} ),\tilde{\sigma}^2,y_{m}) \left( Q  \mathbf{1} {\{y_{m} \in \mathcal{D}\}} + \mathbf{1} {\{y_{m} \not \in \mathcal{D}\}} \right),
\end{eqnarray*}
if $m = n $,
\begin{eqnarray*}
f_{n,m}(x,y) & = & \prod_{i=1}^{n-1} \left(   \phi(y_{i-1} + ( x_i - x_{i-1}),\tilde{\sigma}^2,y_i) (1 - Q)  \mathbf{1} {\{y_i \in \mathcal{D}\}} \right) \\
& & \phi(y_{n-1} + ( x_n - x_{n-1} ),\tilde{\sigma}^2,y_{n}) \\
& & \left( \mathbf{1} {\{y_{n} \not \in \mathcal{D}\}} + (1 - Q) \mathbf{1} {\{y_{n} \in \mathcal{D}\}} \mathbf{1} {\{x_{n} \in \mathcal{D}\}}  +
P \mathbf{1} {\{y_{n} \in \mathcal{D}\}} \mathbf{1} {\{x_{n} \not \in \mathcal{D}\}}  \right),
\end{eqnarray*}
and if $m \geq n+1 $,
\begin{eqnarray*}
f_{n,m}(x,y) & = & \prod_{i=1}^{n-1} \left(    \phi(y_{i-1} + ( x_i - x_{i-1}),\tilde{\sigma}^2,y_i) (1 - Q) \mathbf{1} {\{y_i \in \mathcal{D}\}} \right) \\
& &  \phi(y_{n-1} + ( x_n - x_{n-1} ),\tilde{\sigma}^2,y_{n}) \mathbf{1} {\{y_{n} \in \mathcal{D}\}} \left(    Q \mathbf{1} {\{x_{n} \in \mathcal{D}\}}  +  (1 - P) \mathbf{1} {\{x_{n} \not \in \mathcal{D}\}}   \right) \\
& & \prod_{i=n+1}^{m-1} \left(  \phi(y_{i-1} ,\sigma^2,y_i) (1 - P)  \mathbf{1} {\{y_i \in \mathcal{D}\}}  \right) \\
& & \phi(y_{m-1} ,\sigma^2,y_{m}) \left( \mathbf{1} {\{y_{m} \not \in \mathcal{D}\}}
 + P  \mathbf{1} {\{y_{m} \in \mathcal{D}\}} \right). \\
\end{eqnarray*}

\end{prpstn}

\subsection{Description of the two-dimensional case and expression of the probability density functions} \label{subsection:two:dim:case}

\subsubsection{Description of the neutron transport problem} \label{subsubsection:presentation:two:dim}

The monokinetic particle evolves in $\xR^2$, and its birth takes place at the source point
$s = (-s_x,0)$, with $s_x > 0$.
The domain of interest is a box $B = [-\frac{L}{2},\frac{L}{2}]^2$ with $s_x < \frac{L}{2}$, in which there is an obstacle sphere $S = \left\{x \in \xR^2 ; |x| \leq l \right\}$, with $l<\frac{L}{2}$ and where $|x|$ is the Euclidean norm of $x \in \xR^2$.

We consider two media. The obstacle sphere is composed of ``poison'' and the rest of $\xR^2$ is composed of ``water''. Furthermore if the monokinetic particle leaves the box, it is considered to have gone too far away, and subsequently it is absorbed at the first collision point in the exterior of the box.
The probability of absorption $P(x_n)$ in \eqref{eq:homogeneous:transition:monokinetic} is hence $P(x_n) = \mathbf{1} {\{x_n \not \in B\}} + P_w \mathbf{1} {\{x_n \in B \backslash S\}} +P_p \mathbf{1} {\{x_n \in S\}}$, where $0 \leq P_w \leq 1$ and $0 \leq P_p \leq 1$ are the probabilities of absorption in the water and poison media.

We consider a detector, defined as the sphere $\left\{x \in \xR^2 ; |x - (d_x,0)| \leq l_d \right\}$, with $l < d_x - l_d$ and $d_x + l_d < L/2$, so that the detector is in $B \backslash S$. The event of interest is that the monokinetic particle makes a collision in the detector, before being absorbed.
With $(x_i)_{i \in \xN^*}$ a trajectory of the monokinetic particle and when using the interacting-particle method of Section \ref{section:reminder:AG}, the event of interest is expressed by $\Phi(x) \geq 0$, with $\Phi(x) = l_d - \inf_{ i \in \xN^*; x_i \neq \Delta} |x_i - (d_x,0)|$.
[Note that the probability of absorption in the detector is $P_w$ but that this probability could actually be defined arbitrarily, since it has no impact on the event of interest ``the monokinetic particle makes a collision in the detector''.]

Finally, let us discuss the distribution of the jumps between collision points, corresponding to $q(X_n,X_{n+1})$ in \eqref{eq:homogeneous:transition:monokinetic}. After a scattering, or birth, at $X_n$, of the monokinetic particle, the direction toward which it travels has isotropic distribution. This direction is here denoted $u$, with $u$ a unit two-dimensional vector.
Then, the sampling of the distance to the next collision point $X_{n+1}$ is as follows: First, the distance $\tau$ is sampled from an exponential distribution with rate $\lambda_{w} >0$, if $X_n$ is in the medium ``water'', or $\lambda_p > \lambda_w$ if $X_n$ is in the medium ``poison''. Then, two cases are possible. First, if the sampled distance is so that the monokinetic particle stays in the same medium while it travels this distance, then the next collision point is $X_{n+1} = X_n + \tau u$. Second, if between $X_n$ and $X_n + \tau u$, there is a change of medium, then the monokinetic particle is virtually stopped at the first medium-change point
between $X_n$ and $X_n + \tau u$. At this point, the travel direction remains the same, but the remaining distance to travel is resampled, from the exponential distribution with the rate corresponding to the new medium. These resampling are iterated each time a sampled distance causes a medium-change. The new collision point $X_{n+1}$ is the point reached by the first sampled distance that does not cause a medium change. Note that, in this precise setting with two media, the maximum number of distance sampling between two collision points is three. This can happen in the case where the collision point $X_n$ is in the box but not in the obstacle sphere, where the sampled direction points toward the obstacle sphere, and where toward this direction, the monokinetic particle enters and leaves the obstacle sphere.

The following algorithm, when tuned with source point $(-s_x,0)$, sums up how trajectories can be sampled according to the above description.

\begin{lgrthm} \label{algo:sampling:f:two:dim}

Objective: from a source point $s \in \mathbb{R}^2$ and from the parameters $B$, $S$, $\lambda_{w}$, $\lambda_p$, $P_w$ and $P_p$, sample a trajectory $x$ as described above.

~

\begin{itemize}
\item Set $i=0$, $x_i=s$ and ``state = active''
\item While ``state = active'' do
\begin{itemize}
\item Set $x_{i+1} = x_i$, $\lambda = \lambda_w \mathbf{1} {\{x_i \in \mathbb{R}^2 \backslash S\}} + \lambda_p \mathbf{1} {\{x_i \in  S\}}$ and ``crossing = true''.
\item Sample a vector $v$ from the uniform distribution on the unit sphere of $\mathbb{R}^2$. 
\item While ``crossing = true'' do
\begin{itemize}
\item Sample $r$ from an exponential distribution with rate $\lambda$.
\item If the medium is the same on all the segment $[x_{i+1},x_{i+1} + rv]$, set ``crossing = false'' and $x_{i+1} = x_{i+1} + rv$.
\item Else, set $x_{i+1}$ as the first medium change point when going from $x_{i+1}$ to $x_{i+1} + rv$ on the segment $[x_{i+1},x_{i+1} + rv]$. Set $\lambda = \mathbf{1} {\{\lambda = \lambda_{w}\}} \lambda_p + \mathbf{1} {\{\lambda = \lambda_{p}\}} \lambda_w$.
\end{itemize}
\item If $x_{i+1} \not \in B$
\begin{itemize}
\item Set ``state = inactive''.
\end{itemize}
\item If $x_{i+1 } \in B$
\begin{itemize}
\item With probability $\mathbf{1} {\{x_{i+1} \not \in S\}} P_w + \mathbf{1} {\{x_{i+1} \in S\}} P_p$, set ``state = inactive''.
\end{itemize}
\item Set $i=i+1$.
\end{itemize}
\item Return the infinite sequence $(x_1,...,x_{i},\Delta,...)$.
\end{itemize}

\end{lgrthm}

The pdf corresponding to Algorithm \ref{algo:sampling:f:two:dim}, of a collision point $X_{n+1}$, conditionally to a collision point $X_n$, is given in Proposition \ref{prop:pdf:Xn:Xnp1:2d} below.

Finally, note that the setting described does constitute a simplified but realistic model for a shielding system in neutron transport. Indeed, first exponentially distributed distances (with possible resample after medium change) and uniform directions between collision points correspond to a very classical approximation in neutron transport theory (see e.g. \cite{zoatwo}). 
Second, it is very common to consider simple schemes of the form source-obstacle-detector to evaluate shielding components, either numerically or experimentally \cite{hossain10study}. In our case, the ``water'' medium constitutes a mild obstacle and the ``poison'' medium an important one (larger collision rate and absorption probability). Of course, not all aspects of neutron transport theory, nor exhaustive representations of industrial shielding systems, are tackled here.


\subsubsection{Expression of the probability density function of a trajectory}

We first set some notations for $v,w \in  \xR^2$. We write $[v,w]$ for the segment between $v$ and $w$.
When $v$ is strictly in the interior of $S$ ($|v| < l$) and $w$ is strictly in the exterior of $S$ ($|w| > l$), we let $c(v,w)$ be the unique point in the boundary of $S$ that belongs to $[v,w]$.
Similarly, for $v,w \in \xR^2 \backslash S$ and when $[v,w]$ has a non-empty intersection with $S$, we denote by $c_{1}(v,w)$ and $c_{2}(v,w)$ the two intersection points between $[v,w]$
and the boundary of $S$. The indexes $1$ and $2$ are so that $|v-c_1(v,w)| \leq |v-c_2(v,w)|$.
For $v,w \in \xR^2 \backslash S$, we let $I(v,w)$ be equal to $1$ if $[v,w]$ has a non-empty intersection with $S$ and $0$ otherwise.

The computation of $c(v,w)$, $I(v,w)$, $c_1(v,w)$ and $c_2(v,w)$ are equally needed for a monokinetic-particle simulation (Algorithm \ref{algo:sampling:f:two:dim}), and for the computation of the corresponding pdf of Proposition \ref{prop:pdf:initial:two:dimensional}. The four quantities can be computed explicitly.
We now give the pdf of the collision point $X_{n+1}$, conditionally to a scattering or a birth point $X_n$.

\begin{prpstn} \label{prop:pdf:Xn:Xnp1:2d}
Consider a scattering, or birth, point $x_n \in B$. Then, the pdf of the collision point $X_{n+1}$, conditionally to $x_n$, is denoted $q(x_n,x_{n+1})$ and is given by, if $x_n \in B \backslash S$
\begin{eqnarray*}
  q(x_n,x_{n+1})  &=&  \frac{1}{2 \pi |x_n - x_{n+1}|}   
 \lambda_{w}  \exp{ \left(  - \lambda_{w} |x_n - x_{n+1}| \right) } (1-I(x_n,x_{n+1})) \mathbf{1} {\{x_{n+1} \in \xR^2 \backslash S\}}   \nonumber \\
 && + \frac{1}{2 \pi |x_n - x_{n+1}|}  \exp{ \left( - \lambda_{w} | x_n - c_{1}(x_n,x_{n+1}) |  \right) } \exp{ \left( - \lambda_{p} |c_{1}(x_n,x_{n+1}) - c_{2}(x_n,x_{n+1}) | \right) } \nonumber \\
 && ~ ~  \lambda_{w}  \exp{ \left( - \lambda_{w} |  c_{2}(x_n,x_{n+1}) - x_{n+1} | \right) }  I(x_n,x_{n+1}) \mathbf{1} {\{x_{n+1} \in \xR^2 \backslash S\}} \nonumber \\
 && + \frac{1}{2 \pi |x_n - x_{n+1}|}  
 \exp{ \left( - \lambda_{w} |x_n- c(x_n,x_{n+1})| \right) }
\lambda_{p}  \exp{ \left( - \lambda_{p} |c(x_n,x_{n+1}) - x_{n+1} | \right) } \mathbf{1} {\{x_{n+1} \in  S\}}
\end{eqnarray*}
and, if $x_n \in S$,
\begin{eqnarray*} 
q(x_n,x_{n+1}) & = &   \frac{1}{2 \pi |x_n - x_{n+1}|} 
 \lambda_{p}  \exp{ \left( - \lambda_{p} |x_n- x_{n+1}| \right) }  \mathbf{1} {\{x_{n+1} \in  S\}} \\
&& +  \frac{1}{2 \pi |x_n - x_{n+1}|} 
 \exp{ \left( - \lambda_{p} |x_n- c(x_n,x_{n+1})| \right) }
\lambda_{w}  \exp{ \left( - \lambda_{w} |c(x_n,x_{n+1}) - x_{n+1} | \right) } \mathbf{1} {\{x_{n+1} \in  \xR^2 \backslash S\}}.   \nonumber
\end{eqnarray*}
\end{prpstn}
\begin{proof}
The proposition is obtained by using the properties of the exponential distribution, the definitions of $c(x_n,x_{n+1})$, $I(x_n,x_{n+1})$, $c_1(x_n,x_{n+1})$, and $c_2(x_n,x_{n+1})$ and a two-dimensional polar change of variables. The proof is straightforward
but burdensome.
\end{proof}

Using Proposition \ref{prop:pdf:Xn:Xnp1:2d}, we now give the pdf of the monokinetic-particle trajectories obtained from Algorithm \ref{algo:sampling:f:two:dim}.

\begin{prpstn} \label{prop:pdf:initial:two:dimensional}
The pdf, with respect to $(c_0,\mathcal{S},\Pi)$ of Definition \ref{def:tribu:markov:chain} and Proposition \ref{prop:measure_on_S}, of a trajectory $(x_n)_{n \in \xN^*}$, sampled from Algorithm \ref{algo:sampling:f:two:dim}, is $f(x) = \sum_{n \in \xN^*} \mathbf{1}_{A_{n}}(x) f_n(x)$, with
\begin{eqnarray*}
 f_n(x) & = & \prod_{i=1}^{n-1} \left( q(x_{i-1},x_i) \left[  (1-P_w) \mathbf{1} {\{x_i \in B \backslash S\}} +  (1-P_p)  \mathbf{1} {\{x_i \in  S\}} \right] \right) \\
& & q(x_{n-1},x_n) \left[ \mathbf{1} {\{x_n \not \in B\}} + P_w \mathbf{1} {\{x_n \in B \backslash S\}}  + P_p  \mathbf{1} {\{x_n \in  S\}} \right] ,
\end{eqnarray*}
where $x_0=(-s_x,0)$ by convention, and with $q(x_{i-1},x_i)$ and $q(x_{n-1},x_n)$ as in Proposition \ref{prop:pdf:Xn:Xnp1:2d}.
\end{prpstn}

\subsubsection{Description of the trajectory perturbation method} \label{subsubsection:descriptio:perturbation:2d}

The perturbation method is parameterized by $\tilde{\sigma}^2 >0$, $0 \leq Q_w \leq 1$
and $0 \leq Q_p \leq 1$.
Let us consider a historical trajectory $(x_i)_{i \in \xN^*}$, absorbed at collision $n$.
As in Section \ref{subsection:one:dim:case}, the set of birth and collision points of the perturbed monokinetic particle is an inhomogeneous Markov chain $(Y_i)_{i \in \xN}$, so that $Y_0 = 0$.
We modify independently the collision points of the initial trajectory, and, if the perturbed trajectory outsurvives the initial one, we generate the sequel with the initial distribution. Similarly to Section \ref{subsubsection:one:dim:perturbation:procedure:Pa:non:zero}, we perturb the absorption/non-absorption sampling by changing the initial values with probabilities $Q_w$ and $Q_p$, if the initial and perturbed collision points are both in $B \backslash S$ or both in $S$. If this is not the case, we sample the absorption/non-absorption for the perturbed monokinetic particle with the initial probabilities $P_w$ and $P_p$.

More precisely, for $i \leq n-1$, and if the perturbed monokinetic particle has not been absorbed before collision point $Y_i$, it is absorbed at collision point $Y_i$ with probability $P(x_i,Y_i)$ with

\begin{equation} \label{eq:perturbed:proba:absorption:2d:nm1}
 P(x_i,Y_i) =
  \begin{cases}
   1 & \text{if } Y_{i} \in \xR^2 \backslash B \\
   P_w       & \text{if } Y_i \in B \backslash S \text{ and } x_i \in S \\
      P_p       & \text{if } Y_i \in  S \text{ and } x_i \in B \backslash S \\
            Q_w       & \text{if } Y_i \in  B \backslash S \text{ and } x_i \in B \backslash S \\
            Q_p       & \text{if } Y_i \in   S \text{ and } x_i \in S \\
  \end{cases}.
\end{equation}

Similarly to the one-dimensional case, by taking $Q_w$ smaller than $\min(P_w,1-P_w)$, and $Q_p$ smaller than $\min(P_p,1-P_p)$, we can modify rather mildly the initial trajectories.

If the perturbed monokinetic particle is not absorbed at collision point $Y_i$, its next collision point is $Y_{i+1} = x_{i+1} + \epsilon_{i+1}$, where the
$(\epsilon_i)_{1 \leq i \leq n}$ are independent and where $\epsilon_i$ follows a $\mathcal{N}(0,\tilde{\sigma}^2 I_2)$ distribution, where $I_2$ is the $2 \times 2$ identity matrix.
If the perturbed monokinetic particle has not been absorbed before collision point $Y_n$, then it is absorbed with probability  $P(x_n,Y_n)$ given by

\begin{equation} \label{eq:perturbed:proba:absorption:2d:n}
 P(x_n,Y_n) =
  \begin{cases}
   1 & \text{if } Y_{n} \in \xR^2 \backslash B \\
   P_w       & \text{if } Y_n \in B \backslash S \text{ and } x_n \in S \\
      P_w       & \text{if } Y_n \in B \backslash S \text{ and } x_n \in \xR^2 \backslash B \\
      P_p       & \text{if } Y_n \in  S \text{ and } x_n \in B \backslash S \\
            P_p       & \text{if } Y_n \in  S \text{ and } x_n \in \xR^2 \backslash B \\
        1 - Q_w       & \text{if } Y_n \in B \backslash S \text{ and } x_n \in B \backslash S \\
           1 - Q_p       & \text{if } Y_n \in   S \text{ and } x_n \in S \\
  \end{cases}.
\end{equation}

As in Section \ref{subsection:one:dim:case}, if the perturbed monokinetic particle has not been absorbed before collision point $Y_n$, the sequel of the trajectory of the perturbed monokinetic particle is sampled as the initial particle would be sampled if its collision point $n$ was $Y_n$.

Finally, the following algorithm sums up how perturbed trajectories can be sampled.

\begin{lgrthm} \label{algo:sampling:kappa:two:dim}

Objective: from an initial trajectory $x$ absorbed at collision $n$ and from the parameters $B$, $S$, $\lambda_{w}$, $\lambda_p$, $P_w$, $P_p$, $\tilde{\sigma}^2$, $Q_w$ and $Q_p$, sample a perturbed trajectory $y$ as described above.

~

\begin{itemize}
\item Set $i=0$, $y_i=(-s_x,0)$ and ``state = active''
\item While ``state = active'' and $i +1 < n$ do
\begin{itemize}
\item Sample $y_{i+1}$ from the $\mathcal{N}(x_{i+1},\tilde{\sigma}^2 I_2)$ distribution.
\item With probability $P(x_{i+1},y_{i+1})$ given by \eqref{eq:perturbed:proba:absorption:2d:nm1}, set ``state = inactive''.
\item Set $i=i+1$.
\end{itemize}
\item If ``state = inactive'', stop the algorithm and return the infinite sequence $(y_1,...,y_{i},\Delta,...)$.
\item If ``state = active'' do
\begin{itemize}
\item Sample $y_{i+1}$ from the $\mathcal{N}(x_{i+1} ,\tilde{\sigma}^2 I_2) $ distribution.
\item With probability $P(x_{i+1},y_{i+1})$ given by \eqref{eq:perturbed:proba:absorption:2d:n}, set ``state = inactive''.
\item Set $i=i+1$.
\end{itemize}
\item If ``state = inactive'', stop the algorithm and return the infinite sequence $(y_1,...,y_{i},\Delta,...)$.
\item If ``state = active'' do,
\begin{itemize}
\item  Apply Algorithm \ref{algo:sampling:f:two:dim}, with $s = y_i$ and write $(\tilde{x}_1,...,\tilde{x}_q,\Delta,...)$ for the resulting trajectory.
\item Return the infinite sequence $(y_1,...,y_{i},\tilde{x}_1,...,\tilde{x}_q,\Delta,...)$.
\end{itemize}

\end{itemize}

\end{lgrthm}

\subsubsection{Expression of the probability density function of a perturbed trajectory} \label{subsubsection:two:dim:perturbed:pdf:abs}

\begin{prpstn} \label{prop:pdf:perturbed:two:dimensional}
Let us consider a historical trajectory $(x_i)_{i \in \xN^*}$, absorbed at collision $n$.
Let $y_0 = x_0 = (-s_x,0)$ by convention.
The conditional pdf, with respect to $(c_0,\mathcal{S},\Pi)$ of Definition \ref{def:tribu:markov:chain} and Proposition \ref{prop:measure_on_S}, of a trajectory $(y_n)_{n \in \xN^*}$ sampled from Algorithm \ref{algo:sampling:kappa:two:dim}, is $\kappa(x,y) = \sum_{m \in \xN^*} \mathbf{1}_{A_{m}}(y) f_{n,m}(x,y)$ where the $f_{n,m}$ are given by the following. If $m \leq n-1 $,
\begin{eqnarray*}
f_{n,m}(x,y) & = & \prod_{i=1}^{m-1} \left( \phi(x_{i},\tilde{\sigma}^2 I_2,y_i) \left[ 1 - P(x_i,y_i) \right] \right) \\
& & \phi(x_{m},\tilde{\sigma}^2 I_2,y_m) P(x_m,y_m),
\end{eqnarray*}
with $P(x_i,y_i)$ and $P(x_m,y_m)$ as in \eqref{eq:perturbed:proba:absorption:2d:nm1}. If $m = n $,
\begin{eqnarray*}
f_{n,m}(x,y) & = & \prod_{i=1}^{n-1} \left( \phi(x_{i},\tilde{\sigma}^2 I_2,y_i) \left[ 1 - P(x_i,y_i) \right] \right) \\
& & \phi(x_{n},\tilde{\sigma}^2 I_2,y_n) P(x_n,y_n),
\end{eqnarray*}
with $P(x_i,y_i)$ as in \eqref{eq:perturbed:proba:absorption:2d:nm1}
and $P(x_n,y_n)$ as in \eqref{eq:perturbed:proba:absorption:2d:n}. If $m \geq n+1 $,
\begin{eqnarray*}
f_{n,m}(x,y) & = & \prod_{i=1}^{n-1} \left( \phi(x_{i},\tilde{\sigma}^2 I_2,y_i) \left[ 1 - P(x_i,y_i) \right] \right) \\
& & \phi(x_{n},\tilde{\sigma}^2 I_2,y_n) \left[ 1 - P(x_n,y_n) \right] \\
& &\prod_{i=n+1}^{m-1} \left( q(y_{i-1},y_i) \left[ (1-P_w) \mathbf{1} {\{y_i \in B \backslash S\}}  + (1-P_p) \mathbf{1} {\{y_i \in  S\}}  \right] \right) \\
& &  q(y_{m-1},y_m) \left[ \mathbf{1} {\{y_m \not \in B\}} + P_w \mathbf{1} {\{y_m \in B \backslash S\}}  + P_p \mathbf{1} {\{y_m \in  S\}}  \right],
\end{eqnarray*}
with $P(x_i,y_i)$ as in \eqref{eq:perturbed:proba:absorption:2d:nm1},
$P(x_n,y_n)$ as in \eqref{eq:perturbed:proba:absorption:2d:n} and $q(y_{i-1},y_i)$ and $q(y_{m-1},y_m)$ as in Proposition \ref{prop:pdf:Xn:Xnp1:2d}.

\end{prpstn}

\subsection{Final algorithm for probability estimation} \label{subsection:actual:algo:neutron}

The final algorithm for the one and two-dimensional cases is Algorithm \ref{algo:AG:with:HM}, where the objective functions $\Phi(.)$, the unconditional distributions with pdf $f(.)$, and the instrumental kernels $\kappa(.,.)$ are defined in Sections \ref{subsection:one:dim:case} and \ref{subsection:two:dim:case}.
In order to apply Algorithm \ref{algo:AG:with:HM}, it is hence necessary and sufficient to achieve the five following tasks.

\begin{enumerate}\addtocounter{enumi}{0}
\item Evaluating the objective function $\phi(x)$ for any trajectory $x$.
\item Evaluating the pdf $f(x)$ for any trajectory $x$.
\item Evaluating the conditional pdf $\kappa(x,y)$ for any two trajectories $x$ and $y$.
\item Sampling from the distribution with pdf $f(.)$.
\item Sampling from the distribution with pdf $\kappa(x,.)$, for a fixed trajectory $x$.
\end{enumerate}

In the enumeration above, (1) is straightforward from the expressions of $\phi$ given in Sections \ref{subsection:one:dim:case} and \ref{subsection:two:dim:case}. The task (2) is carried out by using Proposition \ref{prop:pdf:initial:one:dimensional} for the one-dimensional case and Proposition \ref{prop:pdf:initial:two:dimensional} for the two-dimensional case. 
The task (3) is carried out by using Propositions \ref{prop:pdf:perturbed:one:dimensional:Pa:zero} or \ref{prop:pdf:perturbed:one:dimensional:Pa:non:zero} for the one-dimensional case and Proposition \ref{prop:pdf:perturbed:two:dimensional} for the two-dimensional case. The tasks (4) correspond to Algorithm \ref{algo:sampling:f:one:dim} for the one-dimensional case and Algorithm \ref{algo:sampling:f:two:dim} for the two-dimensional case. Finally the tasks (5) correspond to Algorithm \ref{algo:sampling:kappa:one:dim} for the one-dimensional case and Algorithm \ref{algo:sampling:kappa:two:dim} for the two-dimensional case.

\subsection{Proofs for Section \ref{section:one:two:dim:case}}

The proofs are based on the following general Proposition \ref{prop:pdf:inhomogeneous:markov:chains}, giving the expression of pdf for inhomogeneous Markov chains that are absorbed in finite-time.

\begin{prpstn} \label{prop:pdf:inhomogeneous:markov:chains}

Consider a sequence of measurable applications $a_n: \xR^d \to [0,1]$, $n \in \xN$, with $a_0=0$. Consider a sequence $(q_n)_{n \in \xN^*}$ of conditional pdf, that is to say $\forall n, y_{n-1}$, $q_n(y_{n-1},y_n)$ is a pdf on $\xR^d$ with respect to $y_{n}$.
Consider a Markov Chain on a probability space $(\Omega,\mathcal{F},P)$ and with values in $\xR^d \cup \{\Delta\}$, $(Y_n)_{n \in \xN}$, so that $Y_0 = y_o$ a.s, when $y_0$ in a non-zero constant of $\xR^d$. Let, $Y_n$ have the non-homogeneous transition kernel
\begin{equation} \label{eq:inhomogeneous:markov:chain}
k_n(y_{n-1},dy_n) =  \delta_\Delta(dy_n) \mathbf{1} {\{y_{n-1} = \Delta\}} +  \left\{ a_n(y_{n-1}) \delta_\Delta(dy_n) + \left[1 - a_n(y_{n-1})\right] q_n(y_{n-1},y_n) dy_n  \right\} \mathbf{1} {\{y_{n-1} \neq \Delta\}}.
\end{equation}
Assume finally that, almost surely, the Markov Chain $Y_n$ reaches $\Delta$ after a finite time.
Then, the application $\omega \to (Y_i(\omega))_{i \in \xN^*}$ is a random variable on $(c_0,\mathcal{S},\Pi)$ (see Definition \ref{def:tribu:markov:chain} and Proposition \ref{prop:measure_on_S}), with probability density function, for $y=(y_i)_{i \in \xN^*}$,
$f(y) = \sum_{n=1}^{+ \infty} \mathbf{1}_{A_{n}}(y) f_n(y)$, with $A_{n}$ as in \eqref{eq:An:in:Section:four:five} and with
\[
f_n(y) =  \prod_{i=1}^{n} \left[ (1-a_{i-1}(y_{i-1})) q_i(y_{i-1},y_i) \right] a_n(y_n),
\]
where $y_0$ is the constant value of $Y_0$ by convention.
\end{prpstn}
\begin{proof}
Proposition \ref{prop:pdf:inhomogeneous:markov:chains} is proved in the same way as in Section \ref{section:HM:on:Markov:chains}.
\end{proof}

The dynamic \eqref{eq:inhomogeneous:markov:chain} is a time-dependent version of \eqref{eq:homogeneous:transition:monokinetic}. Thus, it addresses the unconditional distribution of the monokinetic particle collision points as well as the conditional one, for the instrumental kernel $\kappa$ (Sections \ref{subsubsection:perturbed:pdf:no:abs} \ref{subsubsection:one:dim:perturbed:pdf:abs}
and \ref{subsubsection:two:dim:perturbed:pdf:abs}).

\begin{proof}[Proof of Proposition \ref{prop:pdf:initial:one:dimensional}]
We apply Proposition \ref{prop:pdf:inhomogeneous:markov:chains} with
\[
a_i(y_{i}) = P \mathbf{1} {\{y_i \in D\}}  + \mathbf{1} {\{y_i \not \in D\}}
\]
and
\[
q_i(y_{i-1},y_i) = \phi(y_{i-1},\sigma^2,y_i).
\]
\end{proof}

\begin{proof}[Proof of Proposition \ref{prop:pdf:perturbed:one:dimensional:Pa:zero}]
We denote $x=(x_i)_{i \in \xN^*}$ the initial trajectory, so that $x \in A_{n}$, and $x_0=0$ by convention.
We apply Proposition \ref{prop:pdf:inhomogeneous:markov:chains} with
\[
a_i(y_{i}) =  \mathbf{1} {\{y_i \not \in D\}},
\]
for $i \geq 1 $,
\[
q_i(y_{i-1},y_i) = \phi(y_{i-1} + x_i - x_{i-1},\tilde{\sigma}^2,y_i),
\]
for $1 \leq i \leq n$
and
\[
q_i(y_{i-1},y_i) = \phi(y_{i-1},\sigma^2,y_i),
\]
for  $ i \geq n+1$.
\end{proof}

\begin{proof}[Proof of Proposition \ref{prop:pdf:perturbed:one:dimensional:Pa:non:zero}]
We denote $x=(x_i)_{i \in \xN^*}$ the initial trajectory, so that $x \in A_{n}$, and $x_0=0$ by convention.
We apply Proposition \ref{prop:pdf:inhomogeneous:markov:chains} with
\[
a_i(y_{i}) = Q \mathbf{1} {\{y_i \in D\}} +  \mathbf{1} {\{ y_i \not \in D \}},
\]
for $1 \leq i \leq n-1$,
\[
a_n(y_{n}) = \mathbf{1} {\{y_i \in D\}} \left( (1-Q) \mathbf{1} {\{x_i \in D\}}  + P \mathbf{1} {\{x_i \not \in D\}}  \right) +  \mathbf{1} {\{y_i \not \in D\}},
\]
\[
a_i(y_{i}) = P \mathbf{1} {\{y_i \in D\}}  +  \mathbf{1} {\{y_i \not \in D\}},
\]
for $i \geq n+1$,
\[
q_i(y_{i-1},y_i) = \phi(y_{i-1} + x_i - x_{i-1},\tilde{\sigma}^2,y_i),
\]
for $1 \leq i \leq n$
and
\[
q_i(y_{i-1},y_i) = \phi(y_{i-1},\sigma^2,y_i),
\]
for  $ i \geq n+1$.
\end{proof}

\begin{proof}[Proof of Proposition \ref{prop:pdf:initial:two:dimensional}]
We apply Proposition \ref{prop:pdf:inhomogeneous:markov:chains} with
\[
a_i(y_{i}) = P_p \mathbf{1} {\{y_i \in S\}}  + P_w \mathbf{1} {\{y_i \in B \backslash S\}}  + \mathbf{1} {\{y_i \in
\xR^2 \backslash B\}}
\]
and
\[
q_i(y_{i-1},y_i) = q(y_{i-1},y_i),
\]
with $q(y_{i-1},y_i)$ as in Proposition \ref{prop:pdf:Xn:Xnp1:2d}.
\end{proof}

\begin{proof}[Proof of Proposition \ref{prop:pdf:perturbed:two:dimensional}]
We denote $x=(x_i)_{i \in \xN^*}$ the initial trajectory, so that $x \in A_{n}$, and $x_0=(-s_x,0)$ by convention.
We apply Proposition \ref{prop:pdf:inhomogeneous:markov:chains} with
\[
a_i(y_{i}) = P(x_i,y_i),
\]
for $1 \leq i \leq n-1$ and with $P(x_i,y_i)$ as in \eqref{eq:perturbed:proba:absorption:2d:nm1},
\[
a_n(y_{n}) = P(x_n,y_n),
\]
with $P(x_n,y_n)$ as in \eqref{eq:perturbed:proba:absorption:2d:n},
\[
a_i(y_{i}) = P_p \mathbf{1} {\{y_i \in S\}}  + P_w \mathbf{1} {\{y_i \in B \backslash S\}}  + \mathbf{1} {\{y_i \in
\xR^2 \backslash B\}},
\]
for $i \geq n+1$
\[
q_i(y_{i-1},y_i) = \phi( x_{i},\tilde{\sigma}^2 I_2,y_i),
\]
for $1 \leq i \leq n$ and
\[
q_i(y_{i-1},y_i) = q(y_{i-1},y_i),
\]
for $ i \geq n+1$ and with $q(y_{i-1},y_i)$ as in Proposition \ref{prop:pdf:Xn:Xnp1:2d}.
\end{proof}

\section{Numerical results in dimension one and two} \label{section:numerical:results}

In this Section \ref{section:numerical:results}, we present numerical results for the interacting-particle method of Section \ref{section:reminder:AG}, in the one and two-dimensional cases of Section \ref{section:one:two:dim:case}. We follow a double objective. First we aim at investigating to what extent the ideal results of the interacting-particle method hold (in term of bias and of theoretical confidence intervals). Second, we want to confirm that, when the objective probability is small, the method outperforms a simple Monte Carlo method.

The simple Monte Carlo method is parameterized by a number of Monte Carlo samples $J$. It consists in generating $J$ independent trajectories $x_1,...,x_{J}$ and in estimating $p$ by the empirical proportion of these trajectories that satisfy the small-probability event.
We denote by $\tilde{p}$ the simple Monte Carlo estimator of $p$.

\subsection{Numerical results in dimension one}

\subsubsection{Features of the interacting-particle method} \label{subsubsection:one:dim:not:small:prob}

We first present a simple one-dimensional setting, with no absorption ($P=0$).
We set for the domain $A = -10$, $B = 1$, and for the variance of the increments
$\sigma^2 =1$. As a result, the probability $p$ to estimate is not small. It is easily estimated to be $p=0.13$ by the simple Monte Carlo method.

For the perturbation method, we set $\tilde{\sigma}^2 = 0.1^2$. This choice may not be optimal, but it is reasonable and can be considered as typical for the implementation of the interacting-particle method in this one-dimensional case.

\begin{figure}
\centering
\begin{tabular}{cc}
\includegraphics[height=7.5cm]{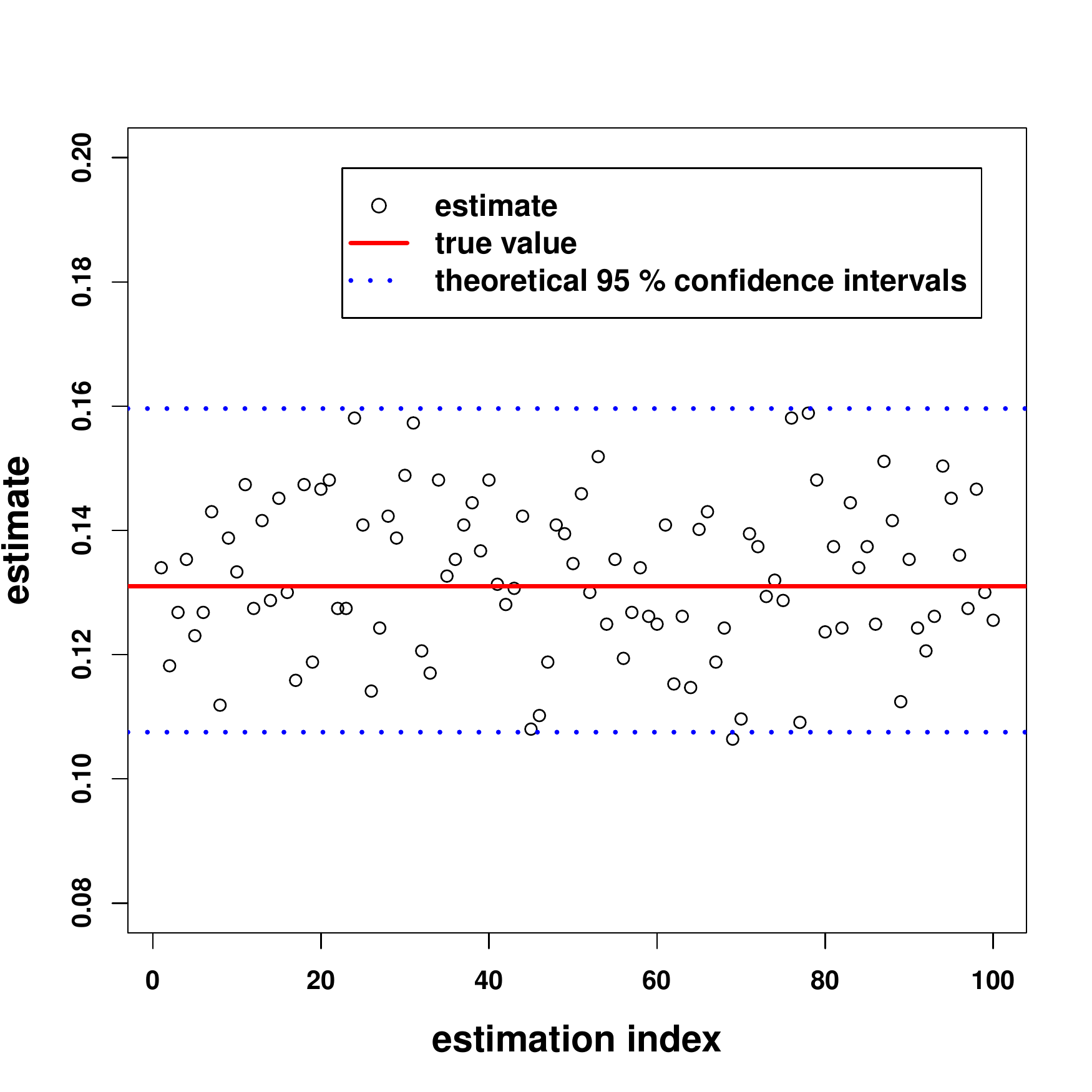}
&
\includegraphics[height=7.5cm]{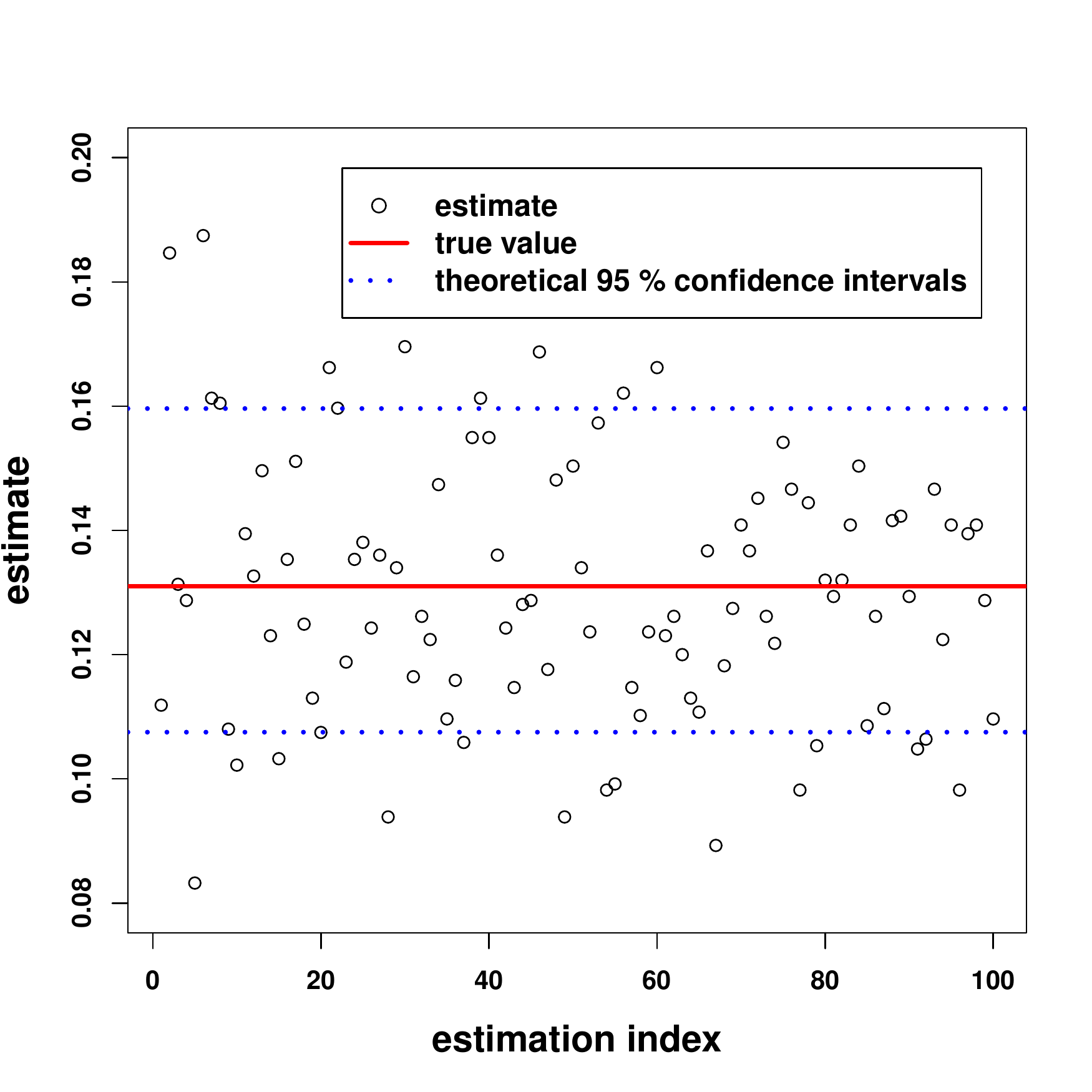}
\end{tabular}
\caption{One-dimensional case. Plot of $100$ independently estimated probabilities with the interacting-particle method \ref{algo:AG:with:HM}, for number of particles $N=200$, and number of iterations in the HM Algorithm \ref{algo:HM} $T=300$ (left) and $T=30$ (right). We also plot the theoretical $95\%$ confidence intervals $I_p$ given by \eqref{eq:IC95:AG} of the case $T=+\infty$. The true probability $p=0.13$ is evaluated quasi-exactly by a simple Monte Carlo method. In both cases, the interacting-particle estimator is empirically unbiased.
For $T=300$, the theoretical confidence interval, obtained in the case $T=+\infty$ is adapted to the practical estimator. For $T=30$ however, the estimator has more variance that the ideal estimator $T=+\infty$ has.}
\label{fig:one:dim:large:proba}
\end{figure}

The results we obtain for $100$ independent estimations for the interacting-particle method are regrouped in Figure \ref{fig:one:dim:large:proba}. We have used $N=200$ particles and $T=300$ and $T=30$ iterations in the HM Algorithm \ref{algo:HM}. Let us first interpret the results for $T=300$ iterations. In this case, we observe that the estimator is empirically non-biased. Furthermore, we also plot the theoretical $95 \%$ confidence intervals for the ideal estimator with $T=+\infty$, that
are approximately (for $N$ large) $
I_p = \left[ p \exp{ \left( -1.96 \sqrt{ \left( \frac{-\log{p}}{N} \right) } \right) } , p \exp{ \left( 1.96 \sqrt{ \left( \frac{-\log{p}}{N} \right) } \right) }\right]$. We also recall from the discussion after \eqref{eq:IC95:AG} that the events $\hat{p} \in I_p$ and $p \in I_{\hat{p}}$ are approximately equivalent when $N$ is large. Hence the coverage probability of $I_p$ for $\hat{p}$ is approximately the probability that $I_{\hat{p}}$ contains $p$, which is the practical quantity of interest. We see on Figure \ref{fig:one:dim:large:proba} that $I_p$ approximately matches the empirical distribution of the estimator $\hat{p}$.
The overall conclusion of this case $T=300$ is that there is a good agreement between theory and practice. This emphasizes the validity of using the interacting-particle method of Algorithm \ref{algo:AG:with:HM}, involving the HM algorithm, in a space that is not a subset of $\xR^d$.

In Figure \ref{fig:one:dim:large:proba}, we also consider the case $T=30$. The estimator is still empirically unbiased. However, its empirical variance is larger, so that the theoretical $95\%$ confidence interval $I_p$ is non-negligibly too thin. This can be interpreted, because when $T$ is small, a new particle at a given conditional sampling step of Algorithm \ref{algo:AG:with:HM} is not independent of the $N-1$ particles that have been kept.
Thus, one can argue that, at each step of Algorithm \ref{algo:AG:with:HM}, the overall set of $N$ particles has more interdependence, so that eventually the estimator has more variance. Nevertheless, on the other hand, an estimation with $T=30$ is $10$ times less time-consuming than an estimation with $T=300$. We further discuss this trade-off problem in Section \ref{subsection:ccl:num:res}.

Finally, for this case of a probability that is not small, we have used simple Monte Carlo as a mean to estimate it quasi-exactly. We have found that the interacting-particle method \ref{algo:AG:with:HM} requires more computation time than the Monte Carlo method, for reaching the same accuracy. We do not elaborate on this fact, since we especially expect the interacting-particle-method to be competitive for estimating a small probability. This is the object of Section \ref{subsubsection:one:dim:small:prob}. For this case of a probability that is not small, we have just investigated the features of the interacting-particle method.

\subsubsection{Comparison with simple Monte Carlo in a small-probability case} \label{subsubsection:one:dim:small:prob}

We now consider a case with possible absorption of the monokinetic particle. Thus we set $P = 0.45$. We keep the same values $\sigma^2 =1$ and $B=1$ as in Section \ref{subsubsection:one:dim:not:small:prob}, but we set $A = -15$. As a result of these parameters for the monokinetic-particle transition kernel, the probability of interest is small. In fact, we have not estimated it with negligible uncertainty. With a simple Monte Carlo estimation of sample size $10^9$, the probability estimate is $ \bar{p} = 6.6 \times 10^{-8}$. We call this estimate the very large Monte Carlo (VLMC) estimate.
Given that the number of successes in this estimate is $66$, which is not very large, we are reluctant to use the Central Limit Theorem approximation for computing $95 \%$ confidence intervals. Instead, we use the Clopper-Pearson interval \cite{clopper34use}, for which the actual coverage probability is always larger than $95 \%$. This $95\%$ confidence interval is there equal to $[5.1 \times 10^{-8},8.4 \times 10^{-8}]$. This uncertainty is small enough for the conclusions we will draw from this case. Finally, note that this very large Monte Carlo estimate is not a benchmark for the interacting-particle method, because it is much more time consuming.

For the interacting-particle method, we set $N=200$ particles, and for the HM algorithm, we set $T=300$ iterations. We use $\tilde{\sigma}^2 = 0.1^2$ and $Q = 0.2$ for the perturbation method. We still denote $\hat{p}$ the obtained estimator for $p$.  We consider a third estimator, that we denote $\tilde{p}$ and that consists in the simple Monte Carlo estimator with sample size $5 \times 10^6$. This sample size is appropriate to compare the efficiency of the interacting-particle and Monte Carlo methods, as we will show below.

The first criterion for comparing the two estimators $\hat{p}$ and $\tilde{p}$
is their computation time. We have two possible ways to make this comparison. First, we can evaluate the complexities of the two methods. The Monte Carlo method requires to perform $5 \times 10^6$ monokinetic-particle simulations. For each proposed perturbation, the interacting-particle method requires to sample one perturbed trajectory, and to compute its unconditional and conditional pdf. This has to be done approximately $T \times \frac{\log{\bar{p}}}{\log{(1-\frac{1}{N})}} \approx 10^6$ times. Thus, from this point of view, the costs of the two methods have the same orders of magnitude. We can not give a more precise comparison, since the trajectories sampled by the two methods do not necessarily have the same length in the mean sense. Furthermore, it is not obvious to compare the computational cost of an initial sampling, with the costs of a conditional sampling and pdf computations.

Hence, we just compare the computational costs of the two methods by considering their actual computational times, for the implementation we have used. Averaged over all the estimations, the time for the interacting-particle method is $58 \%$ of the time for the Monte Carlo method. Hence, we confirm that the computational costs are of the same order of magnitude, the comparison being nevertheless beneficial to the interacting-particle method.

\begin{figure}
\centering
\begin{tabular}{cc}
\includegraphics[height=7.5cm]{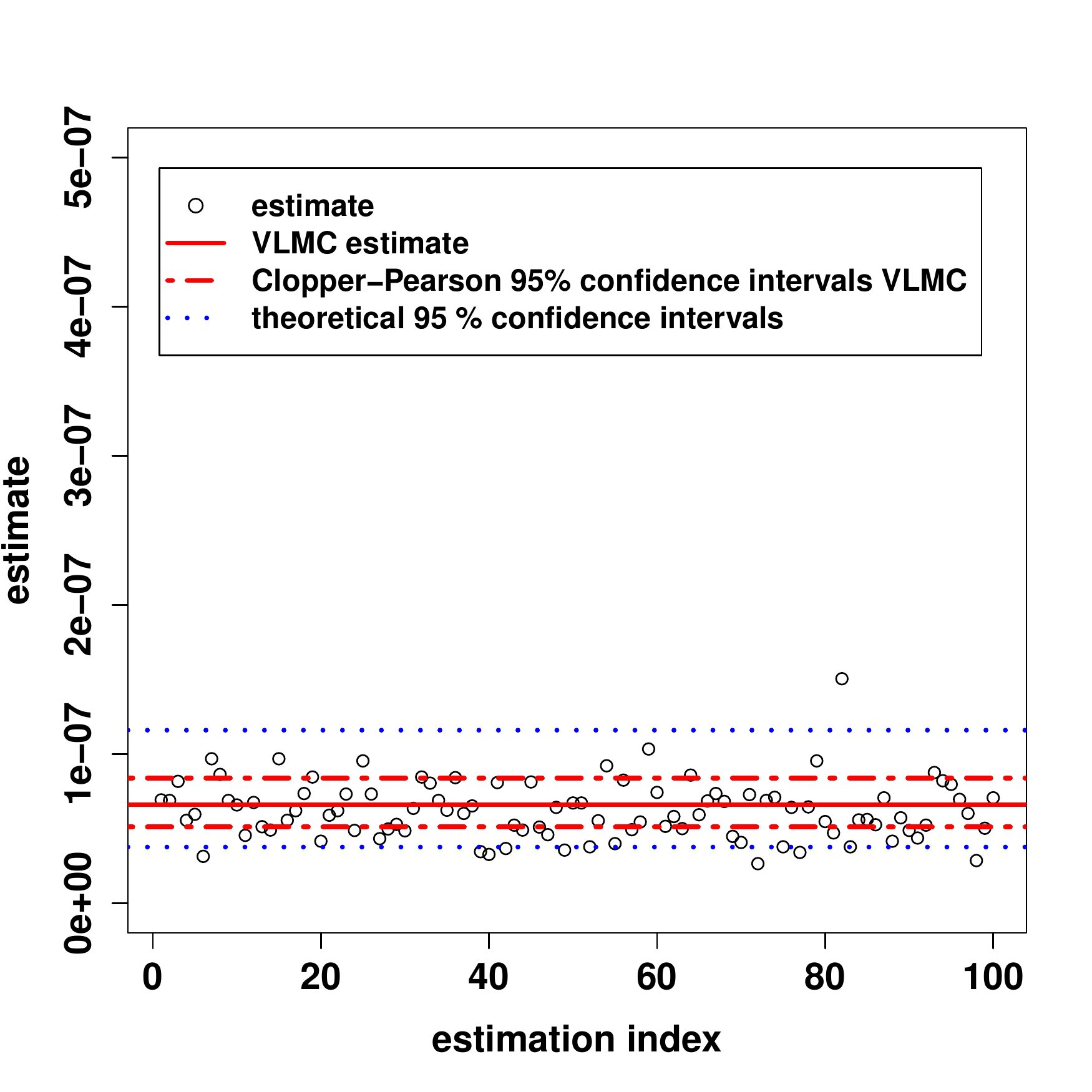}
&
\includegraphics[height=7.5cm]{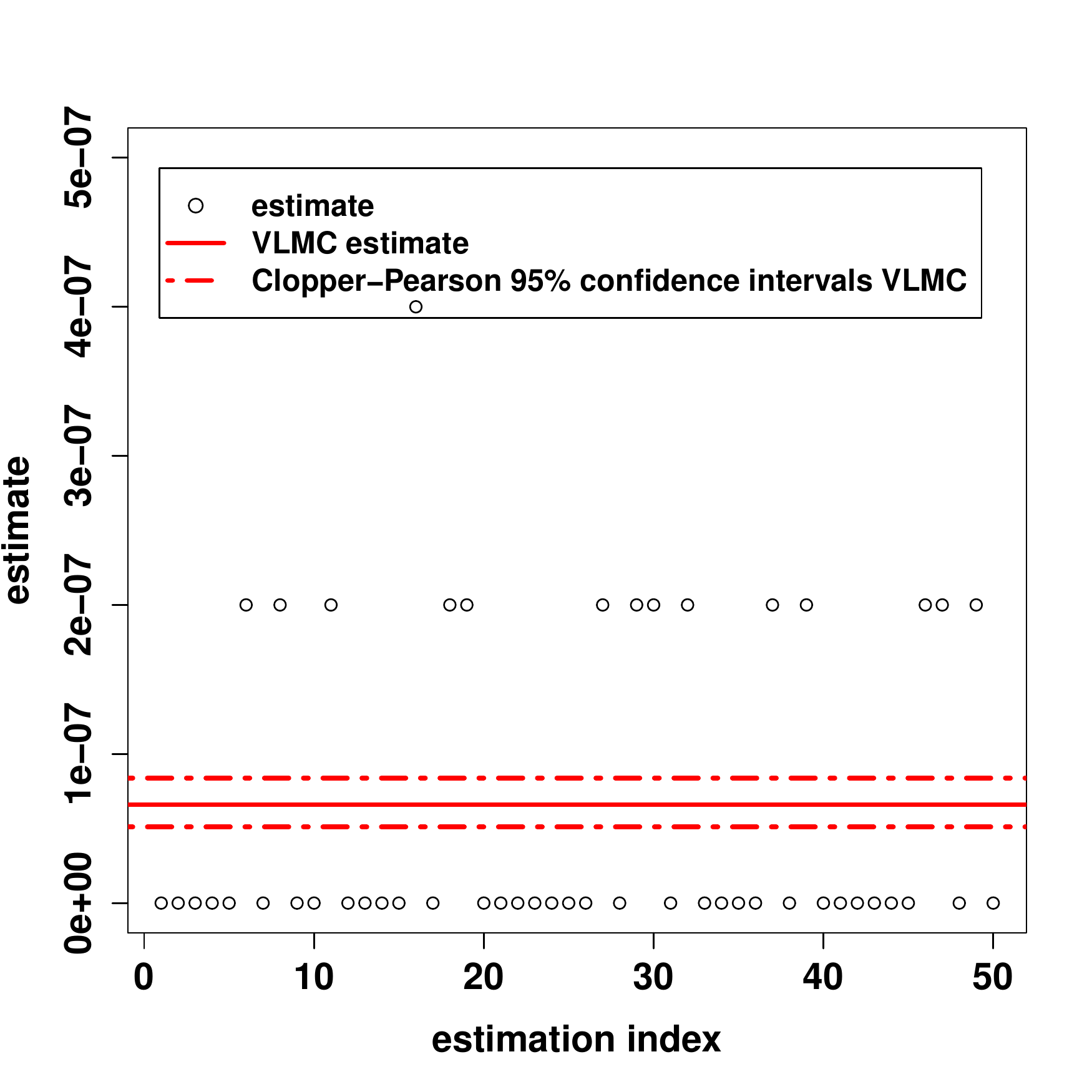}
\end{tabular}
\caption{One-dimensional case. Plot of $100$ independently estimated probabilities with the interacting-particle method \ref{algo:AG:with:HM} for $N=200$ and $T=300$ (left) and $50$ independently estimated probabilities with the Monte Carlo method with sample size $5 \times 10^6$ (right). We plot a very large Monte Carlo estimate of the true probability $\bar{p}=6.6 \times 10^{-8}$, together with the associated Clopper-Pearson $95 \%$ confidence intervals. For the interacting-particle method, we also plot the theoretical $95 \%$ confidence intervals of the case $T=+\infty$, assuming the true probability is the VLMC estimate. The uncertainty on the VLMC estimate of the true value $p$ of the probability is small enough for our conclusions to hold. These conclusions are that the interacting-particle method outperforms the Monte Carlo method (with sample size $5 \times 10^6$), both in term of computation time and of accuracy.}
\label{fig:one:dim:small:proba}
\end{figure}

We now compare the accuracy of the two methods for estimating the true probability $p$.
On Figure \ref{fig:one:dim:small:proba}, we plot the results of $100$ independent estimations for $\hat{p}$ and $50$ independent estimations for $\tilde{p}$. It appears clearly that the interacting-particle method is more precise in this small probability case.
Especially, consider the empirical Root Mean Square Error criterion, for $n$ independent estimates $\check{p}^1,...,\check{p}^n$, for any estimator $\check{p}$ of $p$: $RMSE_{\check{p}} = \sqrt{  \frac{1}{n}\sum_{i=1}^n ( p - \check{p}^i)^2 }$. Regardless of the value of $p$ in the very large Monte Carlo $95 \%$ confidence interval $[5.1 \times 10^{-8},8.4 \times 10^{-8}]$, the RMSE is smaller for $\hat{p}$ than for  $\tilde{p}$. If we assume $p = \bar{p}$, then the RMSE is $10^{-7}$ for $\tilde{p}$ and $2 \times 10^{-8}$ for $\hat{p}$.

A comparison ratio for $\hat{p}$ and $\tilde{p}$, taking into account both computational time and estimation accuracy (in line with the efficiency in \cite{hammersley65monte}), is the quality ratio defined by $\frac{ \sqrt{\mathop{TIME}_{\tilde{p}}} \times \mathop{RMSE}_{\tilde{p}}}{\sqrt{\mathop{TIME}_{\hat{p}}} \times \mathop{RMSE}_{\hat{p}}}$, where the four notations $\mathop{TIME}_{\tilde{p}}$, $\mathop{TIME}_{\hat{p}}$, $\mathop{RMSE}_{\tilde{p}}$ and $\mathop{RMSE}_{\hat{p}}$ are self-explanatory. This ratio is $6.7$ here. This is interpreted as: if the two estimation methods were set as to require the same computational time, then the interacting-particle method would be $6.7$ times more accurate (in term of RMSE) as the Monte Carlo method.

Note that, if we had done the comparison from the point of view of the relative estimation errors, instead of the absolute errors, it would have been even more beneficial to the interacting-particle method. Indeed, assuming again $p = \bar{p} = 6.6 \times 10^{-8}$ for discussion, the interacting-particle method does a maximum relative error of $250\%$. On the other hand, the Monte Carlo estimator takes only $3$ different values in Figure \ref{fig:one:dim:small:proba}.
When it takes value $\frac{2}{5 \times 10^6}$
it does a relative error of $600 \%$, when it takes value $\frac{1}{5 \times 10^6}$
it does a relative error of $300 \%$ and when it takes value $0$ it does an infinite relative error. Alternatively, we can also say that, in the majority of the cases, the Monte Carlo estimator does not see any realization of the rare event, so that it can provide only an overly-conservative upper-bound for $p$.

\subsection{Numerical results in dimension two}

\subsubsection{Features of the interacting-particle method} \label{subsubsection:two:dim:not:small:prob}

We now present the numerical results for the two-dimensional case. We set the absorption probability in the water medium $P_w = 0.2$, the absorption probability in the poison medium $P_p = 0.5$, the dimension of the box $[-\frac{L}{2},-\frac{L}{2}]$ $L=10$, the radius of the obstacle sphere $l=2$, the radius of the detector $l_d=0.5$. The positions of the detector and the source are given by $d_x = s_x = 3$. We set the rate of collisions in the water medium $\lambda_{w} = 0.2$ and in the poison medium $\lambda_{p} = 2$.
As a result, the probability is $p=2\times 10^{-4}$ and is evaluated quasi-exactly by a Monte Carlo sampling, similarly to Section \ref{subsubsection:one:dim:not:small:prob}.

This value is not very small, so that we do not compare the interacting-particle method with the Monte Carlo method. We just aim at showing that the interacting-particle method is valid in this two-dimensional setting,
which is representative of shielding studies with Monte Carlo codes as discussed in Section \ref{section:one:two:dim:case}. 

For the HM perturbations of Algorithm \ref{algo:HM}, we set the collision-point perturbation variance $\tilde{\sigma}^2 = 0.5^2$, the probability of changing the absorption/non absorption in the obstacle sphere $Q_p = 0.1$
and in the rest of the box $Q_w = 0.05$. As in Section \ref{subsubsection:one:dim:not:small:prob}, these settings are reasonable, but are not tuned as to yield an optimal performance of the interacting-particle method.

\begin{figure}
\centering
\includegraphics[height=9cm]{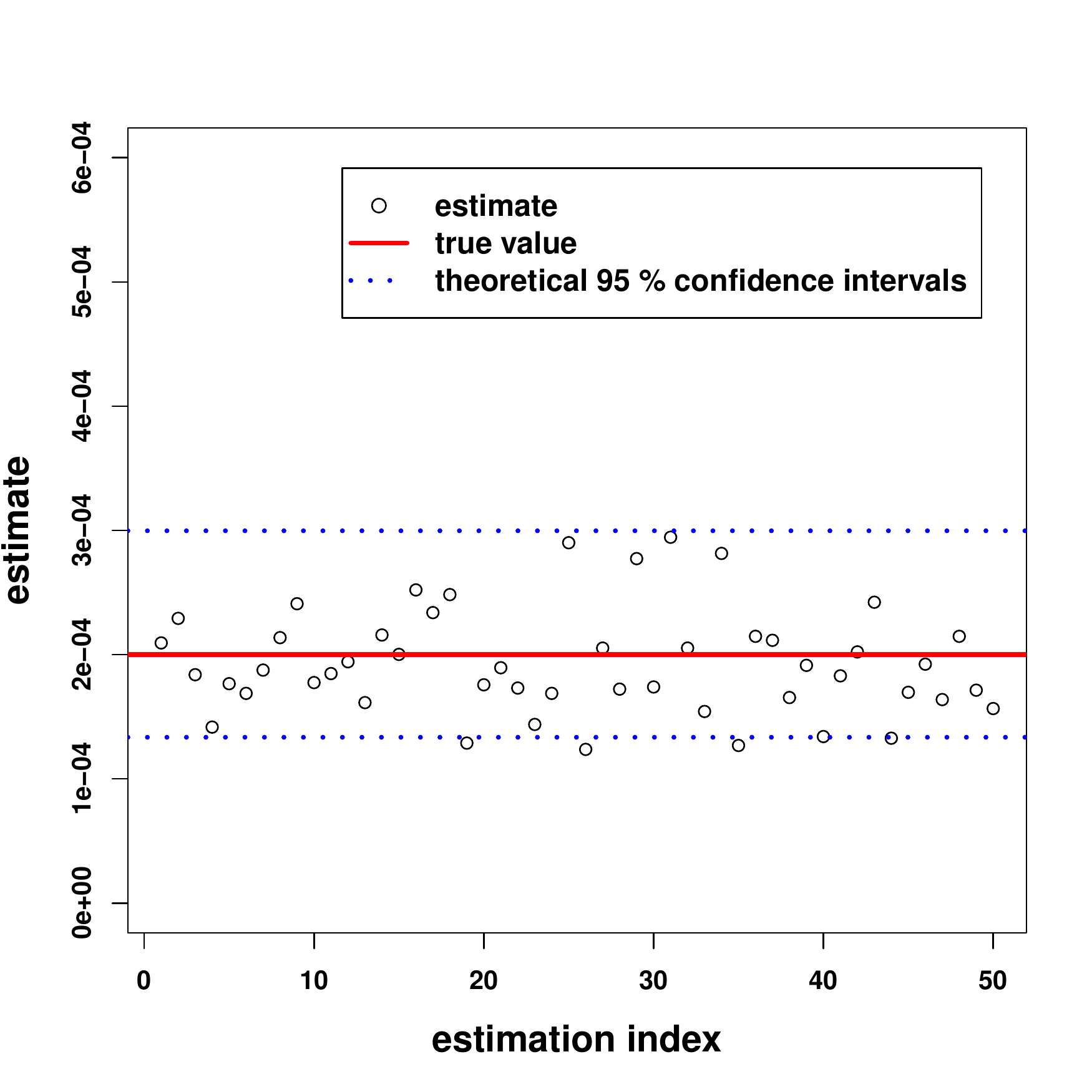}
\caption{Two-dimensional-case. Plot of $50$ independently estimated probabilities with the interacting-particle method \ref{algo:AG:with:HM}, for number of particles $N=200$ and number of iterations in the HM algorithm $T=300$. We also plot the theoretical $95\%$ confidence intervals $I_p$ given by \eqref{eq:IC95:AG} of the case $T=+\infty$. The true probability $p=0.2 \times 10^{-4}$ is evaluated quasi-exactly by a simple Monte Carlo method. The interacting-particle estimator is empirically unbiased and
the $95\%$ theoretical confidence interval, obtained in the case $T=+\infty$, is adapted to the practical estimator.}
\label{fig:two:dim:not:small:proba}
\end{figure}

In Figure \ref{fig:two:dim:not:small:proba}, we present the results for $50$ independent estimations with the interacting-particle method. Empirically, the estimator is unbiased and the theoretical $95 \%$ confidence intervals are valid. This is the same conclusion as in Section \ref{subsubsection:one:dim:not:small:prob}, and is again a validation of the HM algorithm in the space of the monokinetic-particle trajectories.

\subsubsection{Comparison with simple Monte Carlo in a small-probability case} \label{subsubsection:two:dim:small:prob}

We now consider the case of a small probability. For this, we set the absorption probability in the
obstacle sphere $P_p = 0.7$ and in the rest of the box $P_w = 0.5$, the dimension of the box $[-\frac{L}{2},-\frac{L}{2}]$ $L=10$, the radius of the obstacle sphere $l=2.5$, the radius of the detector $l_d=0.5$. The positions of the detector and the source are given by $d_x = s_x = 3$. We set the rate of collisions in the water medium $\lambda_{w} = 2$ and in the poison medium $\lambda_{p} = 3$. In essence, the obstacle sphere is larger than in Section \ref{subsubsection:two:dim:not:small:prob}, the absorption probabilities are larger, and the collision rates are larger, thus yielding all the more frequent absorption.

The probability $p$ is estimated by very large Monte Carlo with sample size $1.25 \times 10^9$. The estimate is $\frac{22}{1.25 \times 10^9} \approx 1.76 \times 10^{-8}$. Similarly to Section \ref{subsubsection:one:dim:small:prob}, the Clopper-Pearson $95\%$ confidence interval for the probability is $[10^{-8},2.5 \times 10^{-8}]$. It is also small enough for validating the discussion that follows.

We compare the estimators $\hat{p}$, with $N=200$ particles and $T=300$ iterations in the HM algorithm, and the estimator $\tilde{p}$ with sample size $5 \times 10^6$. We have found that the computation time for the estimator $\hat{p}$ is, on average, $88 \%$ of that of the estimator $\tilde{p}$.

\begin{figure}
\centering
\begin{tabular}{cc}
\includegraphics[height=7.5cm]{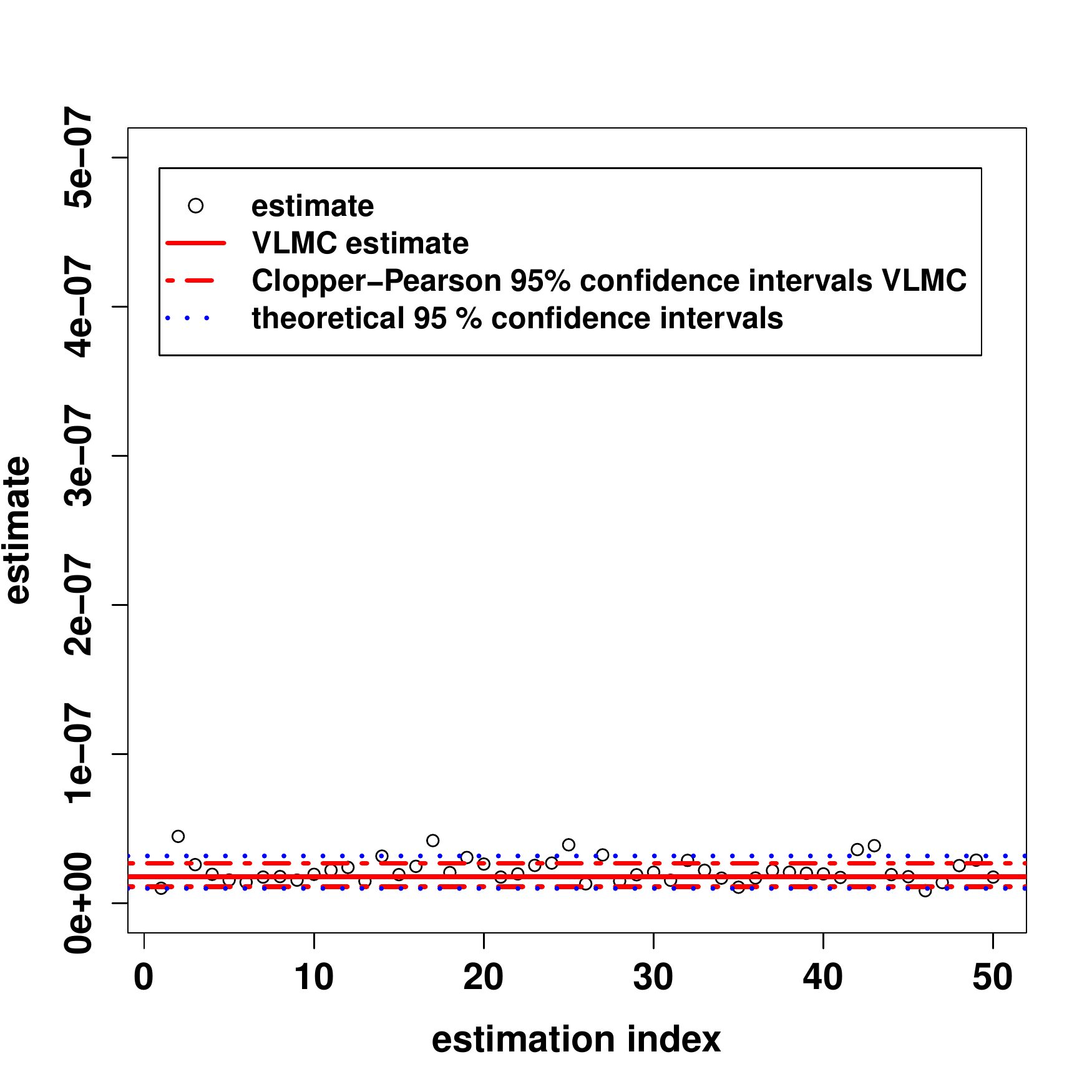}
&
\includegraphics[height=7.5cm]{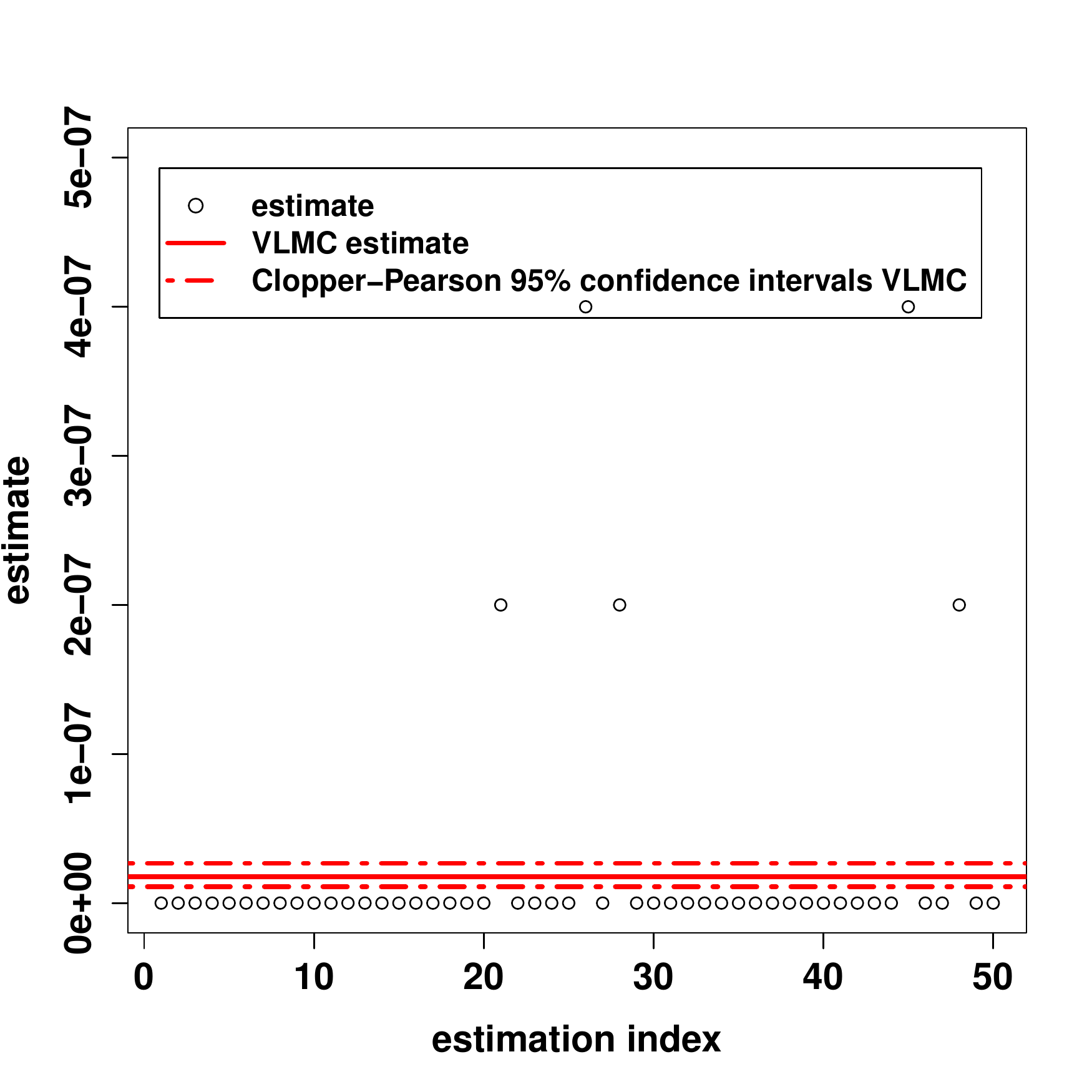}
\end{tabular}
\caption{Two-dimensional-case. Plot of $50$ independently estimated probabilities with the interacting-particle method of Algorithm \ref{algo:AG:with:HM} for $N=200$ and $T=300$ (left) and with the simple Monte Carlo method with sample size $5 \times 10^6$ (right). We plot a very large Monte Carlo estimate of the true probability $\bar{p}=1.76 \times 10^{-8}$, together with the associated Clopper-Pearson $95 \%$ confidence intervals. For the interacting-particle method, we also plot the theoretical $95 \%$ confidence intervals of the case $T=+\infty$, assuming the true probability is the VLMC estimate. The uncertainty on the VLMC estimate of the true value $p$ of the probability is small enough for our conclusions to hold. These conclusions are that the interacting-particle method outperforms the Monte Carlo method (with sample size $5 \times 10^6$), both in term of computation time and of accuracy.}
\label{fig:two:dim:small:proba}
\end{figure}

Now, concerning estimation accuracy, the results are presented in Figure \ref{fig:two:dim:small:proba}. The interacting-particle method outperforms the simple Monte Carlo method, to a greater extent that in Figure \ref{fig:one:dim:small:proba}.
As a confirmation, the quality ratio $\frac{ \sqrt{\mathop{TIME}_{\tilde{p}}} \times \mathop{RMSE}_{\tilde{p}}}{\sqrt{\mathop{TIME}_{\hat{p}}} \times \mathop{RMSE}_{\hat{p}}}$ is $10.5$, against $6.7$ in Figure \ref{fig:one:dim:small:proba}.

\subsection{Discussion on the numerical results} \label{subsection:ccl:num:res}

We now discuss some conclusions on the numerical results of Section \ref{section:numerical:results}. First, in two cases with a probability that is not small (Figures \ref{fig:one:dim:large:proba} and \ref{fig:two:dim:not:small:proba}), the interacting-particle method is empirically unbiased. The theoretical confidence intervals $T=+ \infty$ are in agreement with the empirical distribution for finite $T$, provided that $T$ is large enough.
For the two cases of small probabilities (Figures \ref{fig:one:dim:small:proba} and \ref{fig:two:dim:small:proba}), we do not state conclusions on this question, in one sense or another, because we do not know the probability with negligible uncertainty.

However, for Figures \ref{fig:one:dim:small:proba} and \ref{fig:two:dim:small:proba}, the uncertainty on the probability is small enough to compare the performances of the simple Monte Carlo and interacting-particle methods. The conclusion of this comparison is strongly unilateral, and is that, for a small probability, the interacting-particle method is preferable over a simple Monte Carlo sampling.

We have not carried out numerical test for extremely small probabilities (say, under $10^{-10}$). The reason for that is that we would not have an estimate of these probabilities similar to the very large Monte Carlo estimate $\bar{p}$. That it to say an estimate that comes with confidence intervals with guaranteed coverage probability. Nevertheless, a simple Monte Carlo method, with computational time similar to that of the interacting-particle method, would most likely never see the rare event, and thus only provide an overly conservative upper bound. Thus, the comparison would be even more in favor of the interacting-particle method than for figures
\ref{fig:one:dim:small:proba} and \ref{fig:two:dim:small:proba}.

In Figure \ref{fig:one:dim:large:proba}, we have mentioned the trade-off problem between the number of particles $N$ and the number of HM iterations $T$. The average complexity of the interacting-particle method is proportional to the product $NT$. Naturally, increasing $N$ improves the accuracy of the interacting-particle method. Especially, the variance is proportional to $N$ when $N$ is large, in the ideal case $T=+\infty$. We have seen in Figure \ref{fig:one:dim:small:proba} that increasing $T$ also reduces the variance, which is well interpreted. It is however quite difficult to quantify the dependence between $T$ and the variance of the estimator. We think that the question of this trade-off between $N$ and $T$ would benefit from further investigation.

\begin{figure}
\centering
\includegraphics[height=9cm]{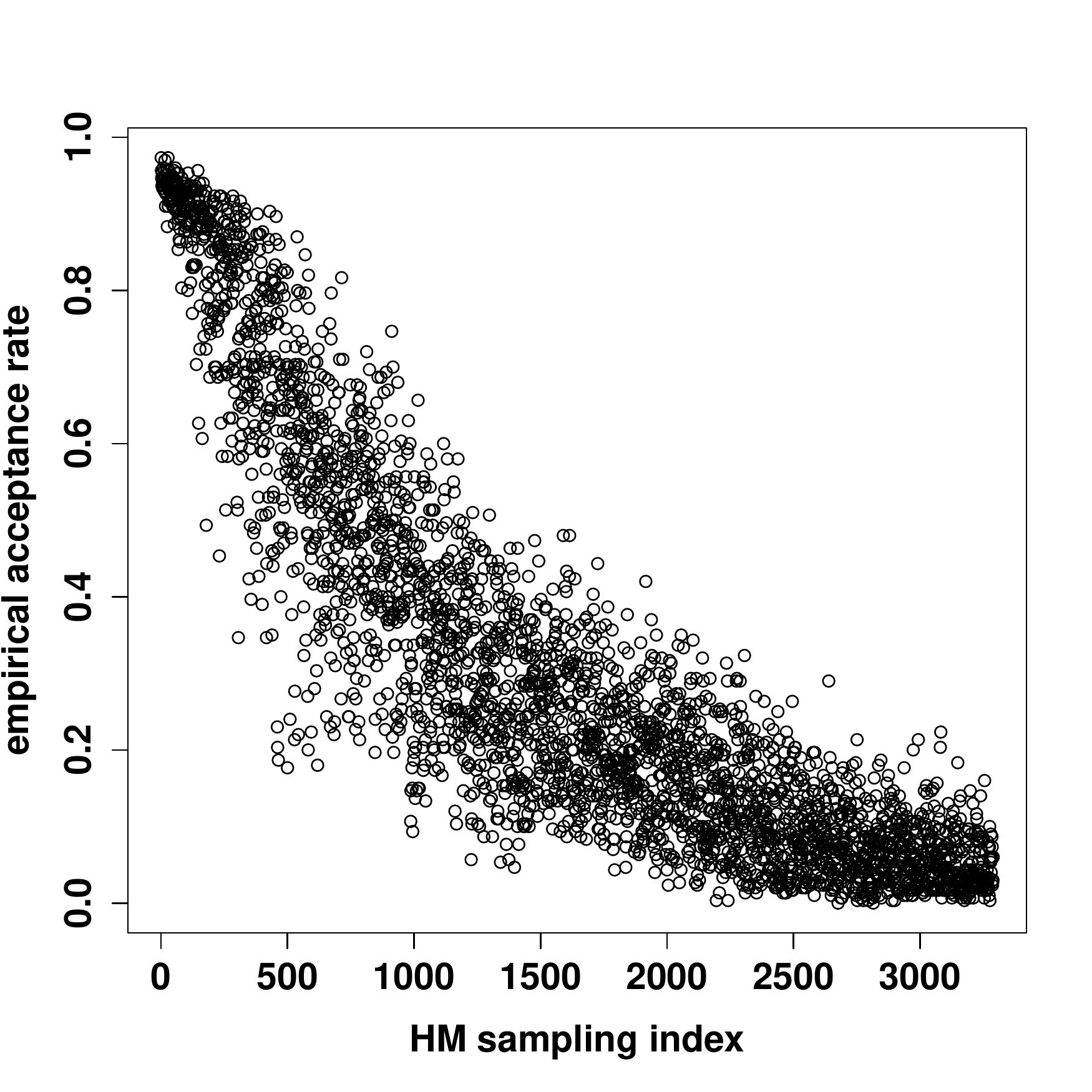}
\caption{Same setting as in Figure \ref{fig:one:dim:small:proba} for the interacting-particle method. We consider one estimation $\hat{p}$ of the interacting-particle method. We plot the empirical acceptance rate, over the $T=300$ proposed perturbations, for each of the $\frac{\log{\hat{p}}}{\log{\left( 1 - \frac{1}{N} \right)}}$ HM samplings of Algorithm \ref{algo:HM}.
The acceptance rate decreases considerably when one gets closer to the rare-event.}
\label{fig:accept:rate}
\end{figure}

In our experiments, we have not optimized the choice of the perturbation method. This would naturally bring a potential additional benefit for the interacting-particle method. Perhaps less natural is the prospect of allowing the perturbation method to vary with the progression of the algorithm. For example, one could use a perturbation method which proposes perturbed trajectories that are closer to the initial ones, when these trajectories are close to the rare event. The results we now present in Figure \ref{fig:accept:rate} support this idea. In Figure \ref{fig:accept:rate}, we plot the acceptance rate in the HM Algorithm \ref{algo:HM} (by acceptance we mean that both the pdf ratio and the objective function conditions are fulfilled), as a function of the progression in the interacting-particle method. This acceptance rate is decreasing, and is small when the interacting-particle method is in the rare-event state. Note that this was not the case in the experiments conducted in \cite{guyader11simulation}.

An other potential tuning of the interacting-particle method is the choice of the objective function $\Phi$, for which the event ``the monokinetic particle makes a collision in the detector'' is equivalent to the event that $\Phi$, evaluated on the trajectory of the monokinetic particle, exceeds a threshold. We have used as a function $\Phi$ the (opposite of the) minimum, over the collision points of the trajectory, of the Euclidean distance to the center of the detector. This choice could be improved. One natural possibility is to replace the Euclidean distance by the optical distance. That is to say the distance traveled in each medium would be weighted by the collision rate in the medium. For some neutron-transport problem, it is also possible to use more specific objective functions, by finding approximations of the importance function, see e.g. \cite{dum}.

\FloatBarrier

\section*{Conclusion}

We have considered the adaptation of the interacting-particle method \cite{guyader11simulation} to a small-probability estimation problem, motivated by shielding studies in neutron transport.
The adaptation is not straightforward, because shielding studies involve working on probability distributions on a set of trajectories that are killed after a finite time.

The contribution brought by the paper it two-fold. First, it has been shown that probability density functions can be defined on this set. This enables to use the Hastings-Metropolis algorithm, which is necessary to implement the method \cite{guyader11simulation} in practice. A convergence result has also been shown for the Hastings-Metropolis algorithm in this setting.

The second contribution of the paper is to give the actual probability density function equations, for implementing the interacting-particle method in an academic one-dimensional problem, and in a simplified but realistic two-dimensional shielding study with monokinetic-particle simulation. In both cases, the method is shown to be valid and to outperform a simple-Monte Carlo estimator, for estimating a small probability.

Prospects are possible for both contributions. First, the proof of the convergence of the Hastings-Metropolis could be extended under more general assumptions. Second, several possibilities for practical improvement of the interacting-particle method are presented in Section \ref{subsection:ccl:num:res}.

\begin{acknowledgement}
The authors are thankful to Josselin Garnier, for suggesting the adaptation of the interacting-particle method to Monte Carlo codes in neutronic and for his advice. They acknowledge the role of Jean-Marc Martinez in the initiation and the conduct of the research project. They thank Fosto Malvagi and Eric Dumonteil for giving them an introduction to neutron-transport problems. Eric Dumonteil also contributed to the writing of the introduction to the neutronic problem. The authors are grateful to Nicolas Champagnat, for his advice during the research project.
They also thank the anonymous reviewers for their helpful suggestions.
This research project was funded by CEA.
\end{acknowledgement}


\bibliographystyle{plain}

\end{document}